\documentstyle[11pt]{article}
\def\Bbb R{{\rm \bf R}}
\def\proclaim#1{\vskip2mm{\bf #1}\em}
\def\endproclaim{\em \vskip2mm}

\def\gathered{\begin{array}{c}}
\def\endgathered{\end{array}}
\def\text{\mbox}

\begin{document}

\title {New Development of the Enclosure Method for Inverse Obstacle Scattering}
\author{Masaru Ikehata\footnote{
This paper is the author version of:
Ikehata, M., New development of the enclosure method for inverse obstacle scattering,
to appear as Chapter {\bf 6} in {\it Inverse Problems and Computational Mechanics} (eds. L. Marin,
L. Munteanu, V. Chiroiu), Vol. {\bf 2}, Editura Academiei, Bucharest, Romania.}
\\
Laboratory of Mathematics, Institute of Engineering\\
Hiroshima University, Higashi-Hiroshima 739-8527, JAPAN}
\date{ }
\maketitle
%\begin{abstract}

%\noindent
%AMS: 35R30

%\noindent KEY WORDS: inverse obstacle scattering problems, enclosure method,
%wave equation, probe method, Helmholtz equation, modified Helmholtz equation
%\end{abstract}

\tableofcontents

\section{Introduction}

The Enclosure Method introduced in \cite{IE2, IE3, IE} has become a well-known {\it guiding principle}
 in attacking various inverse obstacle problems \cite{ISA} governed by partial differential equations.
It is simpler than the Probe Method which has been introduced in \cite{IP, IPR, IP2}.

The Enclosure Method aims at obtaining information about the {\it geometry} of unknown discontinuity.
The method consists of three steps listed below:

$\bullet$  choosing a special solution $v$ depending on a large parameter $\tau>0$ and independent of the unknown discontinuity;

$\bullet$  constructing a so-called indicator function of independent variable $\tau$
by using observation data and $v$;

$\bullet$  studying asymptotic behaviour of the indicator function as
$\tau\longrightarrow\infty$.

From the asymptotic behaviour of the indication function we find a domain that {\it encloses} unknown discontinuity.
The Enclosure Method is quite flexible and its realization depends on the choice of $v$ in the first step
and whether the observation data in the second step depend on $v$ or not.

Now we have many applications of this flexible method to various inverse obstacle problems
governed by elliptic partial differential equations or systems.
See \cite{Isugaku} for the systematic explanation of the Enclosure Method from the beginning and also \cite{ICUBO, ILOG, IkIt, IO, IS, SY, SY2} and \cite{WZ}
with references therein for further applications.
For a nonlinear partial differential equation we cite also a recent remarkable work \cite{BKS}.

It was the paper \cite{IHEAT} which opend the door to various possibilities
of the Enclosure Method in the time domain inverse obstacle problems
governed by the heat or wave equations in one-space dimension.
Now we have several papers \cite{IK1, IW, IK2, IFRAME, IkIt2, IK3} in which the range of application of the Enclosure Method
has been extended to inverse obstacle problems governed by parabolic or hyperbolic equations over a finite time interval
in three-space dimensions.

The aim of this chapter is to make a review of the recent results using the Enclosure Method on
inverse obstacle problems governed by the wave equation and the Maxwell system in time domain.
We also describe some of unsolved problems related to further possibility of the Enclosure Method itself.
Those are not mentioned in the expository paper \cite{Isugaku} and survey paper \cite{ISURVEY}.

\section{The enclosure method for inverse obstacle scattering in time domain}

The description of the problem is simple.
Send a wave and observe the reflected wave by an unknown {\it obstacle} (discontinuity).
What information about the obstacle can one extract from the observed
wave? This type of problems have their origin in sonar, radar, nondestructive testing, etc..

Recently, using the Enclosure Method as the guiding principle,
we considered the problem under the constraints:
sending at most finitely many waves at a finite distance from the obstacle;
observing a reflected wave over a {\it finite time interval} and thus at a {\it finite distance} from the obstacle;
sending and observing places are same.
In particular, we use neither the asymptotic behaviour of the wave as time goes to infinity
nor the far field profile.
In this section we present some of recent results from \cite{IW2, IWCOE} and their applications.

Let $D$ be a non-empty bounded open subset of $\Bbb R^3$ with $C^2$-boundary
such that $\Bbb R^3\setminus\overline D$ is connected.
Let $0<T<\infty$.
Let $f\in L^2(\Bbb R^3)$ satisfy $\text{supp}\,f\cap\overline D=\emptyset$.

We denote by $u_f$ the (weak) solution of the following initial boundary value problem for the wave equation:
%$$\left\{
%\begin{array}{c}
%\displaystyle
%\partial_t^2u-\triangle u=0\,\,\text{in}\,(\Bbb R^3\setminus\overline D)\times\,]0,\,T[,\\
%\\
%\displaystyle
%\frac{\partial u}{\partial\nu}-\gamma(x)\partial_tu-\beta(x) u=0\,\,\text{on}\,\partial D\times\,]0,\,T[,
%\\
%\\
%\displaystyle
%u(x,0)=0\,\,\text{in}\,\Bbb R^3\setminus\overline D,\,\,
%\partial_tu(x,0)=f(x)\,\,\text{in}\,\Bbb R^3\setminus\overline D,
%\end{array}
%\right.
%\tag {2.1}
%$$
\begin{eqnarray}\label{eq:2.1}
\left\{
\begin{array}{ll}
\displaystyle
\partial_t^2u-\triangle u=0 & \text{in} \quad (\Bbb R^3\setminus\overline D)\times\,]0,\,T[,
\\[10pt]
\displaystyle
\frac{\partial u}{\partial\nu}-\gamma(x)\partial_tu-\beta(x) u=0 & \text{on} \quad\partial D\times\,]0,\,T[,
\\[10pt]
\displaystyle
u(x,0)=0 & \text{in} \quad \Bbb R^3\setminus\overline D,
\\[10pt]
\partial_tu(x,0)=f(x) & \text{in} \quad\Bbb R^3\setminus\overline D,
\end{array}
\right.
\end{eqnarray}
where $\nu$ denotes the unit {\it outward normal} to $D$ on $\partial D$, $\beta\in L^{\infty}(\partial D)$,
$\gamma\in L^{\infty}(\partial D)$ and $\gamma(x)\ge 0$ a.e. $x\in\partial D$.
We omit the description about the solution class taken from \cite{DL}.  See \cite{IW2} for the description.

%%%%%%%%Begin Comment%%%%%%%%%%%%
%\textbf{COMMENT~(Liviu Marin):~Instead of the above Latex environment, could you please use the {\sf eqnarray} %environment, i.e.}
%\begin{eqnarray} \label{eq:2.1}
%\left\{
%\begin{array}{ll}
%\displaystyle
%\partial_t^2u-\triangle u=0 & \text{in} \quad (\Bbb R^3\setminus\overline D)\times\,]0,\,T[,
%\\[10pt]
%\displaystyle
%\frac{\partial u}{\partial\nu}-\gamma(x)\partial_tu-\beta(x) u=0 & \text{on} \quad \partial D\times\,]0,\,T[,
%\\[10pt]
%\displaystyle
%u(x,0)=0 & \text{in} \quad \Bbb R^3\setminus\overline D,
%\\[10pt]
%\partial_tu(x,0)=f(x) & \text{in} \quad \Bbb R^3\setminus\overline D,
%\end{array}
%\right.
%\end{eqnarray}
%\textbf{together with a label, i.e. see the Latex command {\sf label} used at the end of the first line in the {\sf eqnarray} %environment, which can be used to refer to this system as follows: For example, see system~($\ref{eq:2.1}$).}

%\textbf{Similar observations hold for the other equations/systems that appear throughout the manuscript and have a %label (in your case, a {\sf tag}).}
%%%%%%%%Begin Comment%%%%%%%%%%%%

The role of $\gamma\ge 0$ can be seen from the formal computation
$$\displaystyle
\mbox{$\cal E$}'(t)=-\int_{\partial D}\gamma(x)\vert\partial_tu\vert^2dx\le 0,
$$
where
$$\displaystyle
\mbox{$\cal E$}(t)=\frac{1}{2}\int_{\Bbb R^3\setminus\overline D}(\vert\partial_tu\vert^2+\vert\nabla u\vert^2)dx
+\frac{1}{2}\int_{\partial D}\beta(x)\vert u\vert^2dS,\,\,
t\in[0,\,T].
$$
We think that the distribution of the values of $\gamma$ and $\beta$ on $\partial D$ is a mathematical model
of the state of the surface of the obstacle.

Let $B$ be the open ball centred at $p$ with {\it very small} radius $\eta$
and satisfy $\overline B\cap\overline D=\emptyset$.
We denote by $\chi_B$ the charactersitic function of $B$.

\noindent
{\bf Problem 2.1.}  Generate $u=u_f$ by the initial data $f=\chi_B$ and observe $u$ on $B$ over
time interval $]0,\,T[$.  Extract information about the geometry of $D$, $\gamma$ and $\beta$ from the observed data.

The correspondence $\displaystyle
(D, \gamma,\beta)\longmapsto u\vert_{B\times\,]0,\,T[}$
is {\it nonlinear}.  Therefore, Problem 2.1 becomes a nonlinear problem.

\subsection{Indicator function}

The Enclosure Method in time domain also introduces an indicator function like the Enclosure Method in frequency domain.
It starts with introducing a special solution with a large parameter.

In what follows we always choose $f=\chi_B$ unless otherwise specified and $u_f$ is the solution of
($\ref{eq:2.1}$).

Let $\tau>0$ and $v=v_f(\,\cdot\,,\tau)\in H^1(\Bbb R^3)$ be the solution of
%$$\displaystyle
%(\triangle-\tau^2)v+f=0\,\,\text{in}\,\Bbb R^3.
%\tag {2.2}
%$$
\begin{eqnarray} \label{eq:2.2}
(\triangle-\tau^2)v+f=0\,\,\text{in}\,\Bbb R^3.
\end{eqnarray}
$v_f$ has the expreesion
%$$\displaystyle
%v_f(x,\tau)=\frac{1}{4\pi}\int_B\frac{e^{-\tau\vert x-y\vert}}{\vert x-y\vert}dy.
%\tag {2.3}
%$$
\begin{eqnarray}\label{eq:2.3}
\displaystyle
v_f(x,\tau)=\frac{1}{4\pi}\int_B\frac{e^{-\tau\vert x-y\vert}}{\vert x-y\vert}dy.
\end{eqnarray}
Define
%$$\displaystyle
%w_f(x,\tau)=\int_0^T e^{-\tau t}
%u_f(x,t)dt,\,\,x\in\Bbb R^3\setminus\overline D.
%\tag {2.4}
%$$
\begin{eqnarray}\label{eq:2.4}
\displaystyle
w_f(x,\tau)=\int_0^T e^{-\tau t}
u_f(x,t)dt,\,\,x\in\Bbb R^3\setminus\overline D.
\end{eqnarray}
Using $v_f$ and $w_f$, we define the indicator function of $\tau$:
$$\displaystyle
I_B(\tau)
= \int_B(w_f-v_f)dx.
$$
This indicator function looks different from the one in the previous version
of the Enclosure Method \cite{IE3}.
So someone may have a question: why should it be called the indicator function?

%However, in \cite{IW2}, we have pointed out the asymptotic formula:
%$$\displaystyle
%I_B(\tau)
%=\int_{\partial\Omega}\left(\frac{\partial v_f}{\partial\nu}w_f
%-\frac{\partial w_f}{\partial\nu}v_f\right)\,dS+O(\tau^{-1}e^{-\tau T}),
%$$
%for an arbitrary fixed $T<\infty$ and provided $\overline B\cap\overline\Omega=\emptyset$
%and $\overline D\subset\Omega$.

However, in \cite{IW2}, it is shown that the asymptotic formula
$$\displaystyle
I_B(\tau)
=\int_{\partial\Omega}\left(\frac{\partial v_f}{\partial\nu}w_f
-\frac{\partial w_f}{\partial\nu}v_f\right)\,dS+O(\tau^{-1}e^{-\tau T}),
$$
is valid for an arbitrary fixed $T<\infty$ and bounded domain $\Omega$ with a smooth boundary, such that
$\overline D\subset\Omega$, $\overline B\cap\overline\Omega=\emptyset$
and $\Bbb R^3\setminus\overline \Omega$ is connected.

And from ($\ref{eq:2.4}$) we have the space-time expression
$$\begin{array}{c}
\displaystyle
\int_{\partial\Omega}\left(\frac{\partial v_f}{\partial\nu}w_f
-\frac{\partial w_f}{\partial\nu}v_f\right)\,dS
=\int_{M}
\left(\frac{\partial (e^{-\tau t} v_f)}{\partial\nu}u_f
-\frac{\partial u_f}{\partial\nu}(e^{-\tau t}v_f)\right)\,dSdt,
\end{array}
$$
where $M=\partial\Omega\times\,]0,\,T[$.
Note also that: ($\ref{eq:2.2}$) implies that $\displaystyle e^{-\tau t} v_f$ satisfies the wave equation in a neighbourhood of $\overline D$.

Theorefore, one can say that $I_B(\tau)$ is essentially similar to the indicator function in the previous version
of the Enlosure Method.  It is a space-time version.

\subsection{Qualitative state of the surface, distance and direction}

The following result is the starting point of the Enclosure Method in time domain.

\proclaim{\noindent Theorem 2.1(\cite{IW2}).}  Let $T>2\text{dist}\,(D,B)$.
Let $C$ be a positive constant.

We have: if $\gamma(x)\le 1-C$ a.e. $x\in\partial D$, then there exists $\tau_0>0$ such that $I_B(\tau)>0$ for all $\tau\ge\tau_0$;
if $\gamma(x)\ge 1+C$ a.e. $x\in\partial D$, then there exists $\tau_0>0$ such that $I_B(\tau)<0$ for all $\tau\ge\tau_0$.
Moreover, in both cases, the formula
%$$\displaystyle
%\lim_{\tau\longrightarrow\infty}\frac{1}{\tau}\log\left\vert I_B(\tau)\right\vert=
%-2\text{dist}\,(D,B),
%\tag {2.5}
%$$
\begin{eqnarray}\label{eq:2.5}
\displaystyle
\lim_{\tau\longrightarrow\infty}\frac{1}{\tau}\log\left\vert I_B(\tau)\right\vert=
-2\text{dist}\,(D,B),
\end{eqnarray}
is valid.
\endproclaim

Define $\displaystyle
d_{\partial D}(p)=\inf_{x\in\partial D}\vert x-p\vert$.
We see that knowing $\text{dist}(D,B)$ is equivalent to
knowing $d_{\partial D}(p)$ since $\text{dist}(D,B)=d_{\partial D}(p)-\eta$.
Thus Theorem 2.1 yields the sphere $\vert x-p\vert=d_{\partial D}(p)$ whose
exterior contains $D$ and on which at least one ponit on $\partial D$ exists.
Thus this should be called an enclosing method by using the exterior of a sphere.
And also, roughly speaking, we can know the {\it qualitative state} of the surface of the unknown obstacle, that is
whether $\gamma>>1$ or $\gamma<<1$ by the {\it signature} of indicator function $I_B(\tau)$ for a large $\tau$
as Theorem 2.1 states.

Finally we present a procedure for making a decision around $p$ whether given direction $\omega\in S^2$
the point $p+d_{\partial D}(p)\omega$ belongs to $\partial D$ or not provided
$d_{\partial D}(p)$ is known.

Fix a large $T$ and small $s\in]0,\,d_{\partial D}(p)[$.
Give direction $\omega\in S^2$ choose an open ball $B'$ centred at $p+s\omega$ such that
$B'$ is contained in the open ball centred at $p$ with radius $d_{\partial D}(p)$.

\noindent
{\bf Step 1.}  Generate $u_f$ by the initial data $f=\chi_{B'}$ and observe $u_f$ on $B'$
over time interval $]0,\,T[$.

\noindent
{\bf Step 2.}  Calculate $d_{\partial D}(p+s\omega)$ from the data obtained in Step 1 via ($\ref{eq:2.5}$) in
Theorem 2.1 in the case when $B$ is replaced with $B'$.

We always have $d_{\partial D}(p+s\omega)\ge d_{\partial D}(p)-s$.  Moreover, it holds that:

$\bullet$  if $d_{\partial D}(p+s\omega)=d_{\partial D}(p)-s$, then $p+d_{\partial D}(p)\omega$ is on $\partial D$;

$\bullet$  if $d_{\partial D}(p+s\omega)>d_{\partial D}(p)-s$, then $p+d_{\partial D}(p)\omega$ is not on $\partial D$.

Therefore one can make a decision around $p$ whether $p+d_{\partial D}(p)\omega$ is on $\partial D$ or not.

\subsection{A sketch of the proof of Theorem 2.1}

The proof of Theorem 2.1 consists of three parts as described below.

\proclaim{\noindent Lemma 2.1.}
We have, as $\tau\longrightarrow\infty$
$$\displaystyle
\Vert R\Vert_{L^2(\Bbb R^3\setminus\overline D)}=O(e^{-\tau\text{dist}\,(D,B)}+e^{-\tau T}),
$$
$$\displaystyle
\Vert\nabla R\Vert_{L^2(\Bbb R^3\setminus\overline D)}=O(\tau(e^{-\tau\text{dist}\,(D,B)}+e^{-\tau T}))
$$
and
$$\displaystyle
\Vert R\Vert_{L^2(\partial D)}=O(\tau^{1/2}(e^{-\tau\text{dist}\,(D,B)}+e^{-\tau T})),
$$
where $R=w_f-v_f$.
\endproclaim

A brief outline of the proof of Lemma 2.1 is as follows.
It follows from ($\ref{eq:2.1}$) and ($\ref{eq:2.2}$) that $R$ satisfies
%$$
%\left\{
%\begin{array}{c}
%\displaystyle
%(\triangle-\tau^2)R=e^{-\tau T}F\,\,\text{in}\,\Bbb R^3\setminus\overline D,
%\\
%\\
%\displaystyle
%\frac{\partial R}{\partial\nu}-cR=-\left(\frac{\partial v}{\partial\nu}-cv\right)+e^{-\tau T}G\,\,\text{on}\,\partial D,
%\end{array}
%\right.
%\tag {2.6}
%$$
\begin{eqnarray}\label{eq:2.6}
\left\{
\begin{array}{ll}
\displaystyle
(\triangle-\tau^2)R=e^{-\tau T}F & \text{in}\quad\Bbb R^3\setminus\overline D,
\\[10pt]
\displaystyle
\frac{\partial R}{\partial\nu}-cR=-\left(\frac{\partial v}{\partial\nu}-cv\right)+e^{-\tau T}G & \text{on}\quad\partial D,
\end{array}
\right.
\end{eqnarray}
where
%$$
%\begin{array}{c}
%\displaystyle
%c=c(x,\tau)=\gamma(x)\tau+\beta(x),\\
%\\
%\displaystyle
%F=F(x,\tau)=\partial_tu(x,T)+\tau u(x,T),\\
%\\
%\displaystyle
%G=G(x)=\gamma(x)u(x,T).
%\end{array}
%\tag {2.7}
%$$
\begin{eqnarray}\label{eq:2.7}
\begin{array}{l}
\displaystyle
c=c(x,\tau)=\gamma(x)\tau+\beta(x),
\\[10pt]
\displaystyle
F=F(x,\tau)=\partial_tu(x,T)+\tau u(x,T),
\\[10pt]
\displaystyle
G=G(x)=\gamma(x)u(x,T).
\end{array}
\end{eqnarray}
($\ref{eq:2.6}$) and integration by parts give
%$$\begin{array}{c}
%\displaystyle
%\int_{\Bbb R^3\setminus\overline D}(\vert\nabla R\vert^2+\tau^2\vert R\vert^2+e^{-\tau T}FR)dx\\
%\\
%\displaystyle
%+\int_{\partial D}
%\left\{c\vert R\vert^2-
%\left(\frac{\partial v}{\partial\nu}-cv\right)R+e^{-\tau T}GR\right\}dS
%=0.
%\end{array}
%\tag {2.8}
%$$
\begin{eqnarray}\label{eq:2.8}
\begin{array}{c}
\displaystyle
\int_{\Bbb R^3\setminus\overline D}(\vert\nabla R\vert^2+\tau^2\vert R\vert^2+e^{-\tau T}FR)dx
\\[12pt]
\displaystyle
+\int_{\partial D}
\left\{c\vert R\vert^2-
\left(\frac{\partial v}{\partial\nu}-cv\right)R+e^{-\tau T}GR\right\}dS
=0.
\end{array}
\end{eqnarray}
Using a trace theorem \cite{G} and the assumption $\gamma(x)\ge 0$ a.e. $x\in\partial D$,
from ($\ref{eq:2.8}$) one can easily deduce the conclusion.  See Lemma 2.1 and (2.28) in \cite{IW2}.

Let us continue the sketch of the proof of Theorem 2.1.
Since $v$ satisfies ($\ref{eq:2.2}$), we obtain
%$$\displaystyle
%\int_{\Bbb R^3\setminus\overline D}fRdx
%=\int_{\partial D}\frac{\partial v}{\partial\nu}RdS
%+\int_{\Bbb R^3\setminus\overline D}(\nabla v\cdot\nabla R+\tau^2 vR)dx.
%\tag {2.9}
%$$
\begin{eqnarray}\label{eq:2.9}
\displaystyle
\int_{\Bbb R^3\setminus\overline D}fRdx
=\int_{\partial D}\frac{\partial v}{\partial\nu}RdS
+\int_{\Bbb R^3\setminus\overline D}(\nabla v\cdot\nabla R+\tau^2 vR)dx.
\end{eqnarray}
On the other hand, from ($\ref{eq:2.6}$) we obtain
%$$\begin{array}{c}
%\displaystyle
%0=\int_{\partial D}\left\{cR-\left(\frac{\partial v}{\partial\nu}-cv\right)\right\}vdS
%+\int_{\Bbb R^3\setminus\overline D}(\nabla R\cdot\nabla v+\tau^2Rv)dx\\
%\\
%\displaystyle
%+e^{-\tau T}
%\left(\int_{\Bbb R^3\setminus\overline D}Fvdx+\int_{\partial D}GvdS\right).
%\end{array}
%\tag {2.10}
%$$
\begin{eqnarray}\label{eq:2.10}
\begin{array}{c}
\displaystyle
0=\int_{\partial D}\left\{cR-\left(\frac{\partial v}{\partial\nu}-cv\right)\right\}vdS
+\int_{\Bbb R^3\setminus\overline D}(\nabla R\cdot\nabla v+\tau^2Rv)dx
\\[12pt]
\displaystyle
+e^{-\tau T}
\left(\int_{\Bbb R^3\setminus\overline D}Fvdx+\int_{\partial D}GvdS\right).
\end{array}
\end{eqnarray}
Taking the difference of ($\ref{eq:2.9}$) from ($\ref{eq:2.10}$), we obtain
$$\begin{array}{c}
\displaystyle
I_B(\tau)
=\int_{\partial D}\left(\frac{\partial v}{\partial\nu}-cv\right)RdS+\int_{\partial D}\left(\frac{\partial v}{\partial\nu}-cv\right)vdS
\\[12pt]
\displaystyle
-e^{-\tau T}
\left(\int_{\Bbb R^3\setminus\overline D}Fvdx+\int_{\partial D}GvdS\right).
\end{array}
%\tag {2.11}
$$
Then, applying ($\ref{eq:2.8}$) to the first term on this right-hand side, we obtain another expression
%$$\begin{array}{c}
%$\displaystyle
%I_B(\tau)
%=\int_{\Bbb R^3\setminus\overline D}(\vert\nabla R\vert^2+\tau^2\vert R\vert^2)dx
%+\int_{\partial D}c\vert R\vert^2dS
%+\int_{\partial D}\left(\frac{\partial v}{\partial\nu}-cv\right)vdS\\
%\\
%\displaystyle
%+e^{-\tau T}\left(\int_{\Bbb R^3\setminus\overline D}FRdx+\int_{\partial D}GRdS
%-\int_{\Bbb R^3\setminus\overline D}Fvdx-\int_{\partial D}GvdS\right).
%\end{array}
%\tag {2.12}
%$$
\begin{eqnarray}\label{eq:2.12}
\begin{array}{c}
\displaystyle
I_B(\tau)
=\int_{\Bbb R^3\setminus\overline D}(\vert\nabla R\vert^2+\tau^2\vert R\vert^2)dx
+\int_{\partial D}c\vert R\vert^2dS
+\int_{\partial D}\left(\frac{\partial v}{\partial\nu}-cv\right)vdS
\\[12pt]
\displaystyle
+e^{-\tau T}\left(\int_{\Bbb R^3\setminus\overline D}FRdx+\int_{\partial D}GRdS
-\int_{\Bbb R^3\setminus\overline D}Fvdx-\int_{\partial D}GvdS\right).
\end{array}
\end{eqnarray}
It is easy to see that, as $\tau\longrightarrow\infty$ $\displaystyle
\Vert v\Vert_{L^2(\Bbb R^3\setminus\overline D)}=O(\tau^{-2})$ and
$\Vert v\Vert_{L^2(\partial D)}=O(e^{-\tau\text{dist}\,(D,B)})$.
Thus from this, ($\ref{eq:2.7}$) and Lemma 2.1 we see that the last term in the right-hand side
on ($\ref{eq:2.12}$) has bound $O(\tau^{-1}e^{-\tau T})$ as $\tau\longrightarrow\infty$.
Therefore  we have, as $\tau\longrightarrow\infty$
%$$\begin{array}{c}
%\displaystyle
%I_B(\tau)
%=E(\tau)+J(\tau)
%+O(\tau^{-1}e^{-\tau T}),
%\end{array}
%\tag {2.13}
%$$
\begin{eqnarray}\label{eq:2.13}
\displaystyle
I_B(\tau)
=E(\tau)+J(\tau)
+O(\tau^{-1}e^{-\tau T}),
\end{eqnarray}
where
$$\displaystyle
E(\tau)=
\int_{\Bbb R^3\setminus\overline D}(\vert\nabla R\vert^2+\tau^2\vert R\vert^2)dx
+\int_{\partial D}c\vert R\vert^2dS
$$
and
$$\displaystyle
J(\tau)
=\int_{\partial D}\left(\frac{\partial v}{\partial\nu}-cv\right)vdS.
$$

A combination of this and Lemma 2.1 yields the following estimates on the bound of
the indicator function.

\proclaim{\noindent Lemma 2.2.}  We have the following asymptotic estimates:

(i)  if $0\le \gamma(x)$ a.e. $x\in\partial D$, then as $\tau\longrightarrow\infty$ we have
%$$\begin{array}{c}
%\displaystyle
%I_B(\tau)
%\ge
%J(\tau)+O(\tau^{-1}e^{-\tau T});
%\end{array}
%\tag {2.14}
%$$
\begin{eqnarray}\label{eq:2.14}
\displaystyle
I_B(\tau)
\ge
J(\tau)+O(\tau^{-1}e^{-\tau T});
\end{eqnarray}

(ii)  if $\gamma(x)\ge C'$ a.e. $x\in\partial D$ for a positive constant $C'$, then as $\tau\longrightarrow\infty$ we have
%$$\begin{array}{c}
%\displaystyle
%I_B(\tau)
%\le J(\tau)
%+\int_{\partial D}
%\frac{1}{c}\left\vert\frac{\partial v}{\partial\nu}-cv\right\vert^2dS
%+O(\tau^{-1}e^{-\tau T});
%\end{array}
%\tag {2.15}
%$$
\begin{eqnarray}\label{eq:2.15}
\displaystyle
I_B(\tau)
\le J(\tau)
+\int_{\partial D}
\frac{1}{c}\left\vert\frac{\partial v}{\partial\nu}-cv\right\vert^2dS
+O(\tau^{-1}e^{-\tau T});
\end{eqnarray}

(iii)  as $\tau\longrightarrow\infty$ we have
%$$\displaystyle
%\vert I_B(\tau)\vert=O(\tau^{2}e^{-2\tau\text{dist}\,(D,B))}+\tau^{-1}e^{-\tau T}).
%\tag {2.16}
%$$
\begin{eqnarray}\label{eq:2.16}
\displaystyle
\vert I_B(\tau)\vert=O(\tau^{2}e^{-2\tau\text{dist}\,(D,B))}+\tau^{-1}e^{-\tau T}).
\end{eqnarray}

\endproclaim

A brief outline of the proof of Lemma 2.2 is as follows.
From ($\ref{eq:2.3}$), ($\ref{eq:2.13}$) and Lemma 2.1 we have ($\ref{eq:2.16}$); ($\ref{eq:2.14}$) is clear from ($\ref{eq:2.13}$)
and the positivity of $E(\tau)$ for $\tau>>1$ which is a consequence of the trace theorem
\cite{G}.
We present here a sketch of the proof of ($\ref{eq:2.15}$).  See \cite{IW2} for the full proof.
Assume that $\gamma(x)\ge C'$ a.e. $x\in\partial D$ for a positive constant $C'$.
Rewrite ($\ref{eq:2.8}$) as
%$$\begin{array}{c}
%\displaystyle
%\int_{\Bbb R^3\setminus\overline D}
%\left(2\vert\nabla R\vert^2
%+2\tau^2\left\vert R+\frac{e^{-\tau T}F}{2\tau^2}\right\vert^2\right)dx\\
%\\
%\displaystyle
%+\int_{\partial D}2c\left\vert R-\frac{1}{2c}\left\{\left(\frac{\partial v}{\partial\nu}-cv\right)-e^{-\tau T}G\right\}\right\vert^2dS\\
%\\
%\displaystyle
%=\frac{e^{-2\tau T}}{2\tau^2}\int_{\Bbb R^3\setminus\overline D}\vert F\vert^2dx
%+\int_{\partial D}\frac{1}{2c}\left\vert\left(\frac{\partial v}{\partial\nu}-cv\right)-e^{-\tau T}G\right\vert^2dS.
%\end{array}
%\tag {2.17}
%$$
\begin{eqnarray}\label{eq:2.17}
\begin{array}{c}
\displaystyle
\int_{\Bbb R^3\setminus\overline D}
\left(2\vert\nabla R\vert^2
+2\tau^2\left\vert R+\frac{e^{-\tau T}F}{2\tau^2}\right\vert^2\right)dx
\\
\\
\displaystyle
+\int_{\partial D}2c\left\vert R-\frac{1}{2c}\left\{\left(\frac{\partial v}{\partial\nu}-cv\right)-e^{-\tau T}G\right\}\right\vert^2dS
\\
\\
\displaystyle
=\frac{e^{-2\tau T}}{2\tau^2}\int_{\Bbb R^3\setminus\overline D}\vert F\vert^2dx
+\int_{\partial D}\frac{1}{2c}\left\vert\left(\frac{\partial v}{\partial\nu}-cv\right)-e^{-\tau T}G\right\vert^2dS.
\end{array}
\end{eqnarray}
Since we have
$$\displaystyle
\tau^2\vert R\vert^2\le 2\tau^2\left\vert R+\frac{e^{-\tau T}F}{2\tau^2}\right\vert^2+\frac{e^{-2\tau T}\vert F\vert^2}{2\tau^2}
$$
and
$$\displaystyle
c\vert R\vert^2\le 2c\left\vert R-\frac{1}{2c}\left\{\left(\frac{\partial v}{\partial\nu}-cv\right)-e^{-\tau T}G\right\}\right\vert^2
+\frac{1}{2c}\left\vert\left(\frac{\partial v}{\partial\nu}-cv\right)-e^{-\tau T}G\right\vert^2,
$$
noting the trivial inequality $\vert\nabla R\vert^2\le 2\vert\nabla R\vert^2$,
from ($\ref{eq:2.17}$) we obtain
%$$\begin{array}{c}
%\displaystyle
%E(\tau)
%\le\frac{e^{-2\tau T}}{\tau^2}\int_{\Bbb R^3\setminus\overline D}\vert F\vert^2dx
%+\int_{\partial D}\frac{1}{c}\left\vert\left(\frac{\partial v}{\partial\nu}-cv\right)-e^{-\tau T}G\right\vert^2dS.
%\end{array}
%\tag {2.18}
%$$
\begin{eqnarray}\label{eq:2.18}
\displaystyle
E(\tau)
\le\frac{e^{-2\tau T}}{\tau^2}\int_{\Bbb R^3\setminus\overline D}\vert F\vert^2dx
+\int_{\partial D}\frac{1}{c}\left\vert\left(\frac{\partial v}{\partial\nu}-cv\right)-e^{-\tau T}G\right\vert^2dS.
\end{eqnarray}
Writing
$$\displaystyle
\left\vert\left(\frac{\partial v}{\partial\nu}-cv\right)-e^{-\tau T}G\right\vert^2
=\left\vert\frac{\partial v}{\partial\nu}-cv\right\vert^2
-2\left(\frac{\partial v}{\partial\nu}-cv\right)e^{-\tau T}G
+e^{-2\tau T}\vert G\vert^2,
$$
we have
$$\begin{array}{c}
\displaystyle
\int_{\partial D}\frac{1}{c}\left\vert\left(\frac{\partial v}{\partial\nu}-cv\right)-e^{-\tau T}G\right\vert^2dS
=\int_{\partial D}\frac{1}{c}\left\vert\frac{\partial v}{\partial\nu}-cv\right\vert^2dS\\
\\
\displaystyle
+O(\tau e^{-\tau(\text{dist}\,(D,B)+T)}+e^{-2\tau T}).
\end{array}
$$
Therefore ($\ref{eq:2.18}$) yields
%$$\begin{array}{c}
%\displaystyle
%E(\tau)
%\le
%\int_{\partial D}\frac{1}{c}\left\vert\frac{\partial v}{\partial\nu}-cv\right\vert^2dS
%+O(\tau e^{-\tau(\text{dist}\,(D,B)+T)}+e^{-2\tau T}).
%\end{array}
%\tag {2.19}
%$$
\begin{eqnarray}\label{eq:2.19}
\displaystyle
E(\tau)
\le
\int_{\partial D}\frac{1}{c}\left\vert\frac{\partial v}{\partial\nu}-cv\right\vert^2dS
+O(\tau e^{-\tau(\text{dist}\,(D,B)+T)}+e^{-2\tau T}).
\end{eqnarray}
Now a combination of ($\ref{eq:2.19}$) and ($\ref{eq:2.13}$) gives ($\ref{eq:2.15}$).

By virtue of Lemma 2.2 it suffices to study the asymptotic behaviour of two Laplace type integrals in (i) and (ii) from below and above,
respectively.  For this we have the following estimates.

\proclaim{\noindent Lemma 2.3.}  Let $C$ be a positivie constant.

(i)  If $0\le \gamma(x)\le 1-C$ a.e. $x\in\partial D$, then there exist positive numbers $\mu$,
$C'$ and $\tau_0$ such that, for all $\tau\ge\tau_0$

$$\displaystyle
J(\tau)\ge C'\tau^{-\mu}e^{-2\tau\text{dist}\,(D,B)}.
$$

(ii)  If $\gamma(x)\ge 1+C$ a.e. $x\in\partial D$, then there exist positive numbers $\mu$,
$C'$ and $\tau_0$ such that, for all $\tau\ge\tau_0$
$$\displaystyle
J(\tau)+\int_{\partial D}
\frac{1}{c}\left\vert\frac{\partial v}{\partial\nu}-cv\right\vert^2dS
\le -C'\tau^{-\mu}e^{-2\tau\text{dist}\,(D,B)}.
$$

\endproclaim

For the proof of Lemma 2.3 we refer the reader to \cite{IW2}.  The proof given therein
is based on an argument done in \cite{IK2} and covers more general $f$.
Here we describe roughly why $\gamma(x)=1$ is exceptional in Lemma 2.3.

Applying the mean value theorem \cite{CH} to ($\ref{eq:2.3}$), we have
$$\displaystyle
v(x)=\frac{\varphi(\tau\eta)}{\tau^3}\frac{e^{-\tau\vert x-p\vert}}{\vert x-p\vert},\,x\in\Bbb R^3\setminus\overline B,
$$
where $\varphi(\xi)=\xi\cosh\xi-\sinh\xi$.

Define $\displaystyle \Lambda_{\partial D}(p) = \big\{ q \in \partial D \, \big\vert \, \vert q - p\vert = d_{\partial D}(p) \big\}$. 
We call $\Lambda_{\partial D}(p)$ the first reflector from $p$ to $\partial D$ and 
the points in the first reflector are called the first-reflection points, going from $p$ to $\partial D$.

Let $x\in\Lambda_{\partial D}(p)$.
Since $\nu_x=(p-x)/\vert x-p\vert$, from the expression above we obtain
$$\displaystyle
\frac{\partial v}{\partial\nu}=\tau v\left(1+\frac{1}{\tau\vert x-p\vert}\right)\sim\tau v
$$
and hence
$$\displaystyle
\frac{\partial v}{\partial\nu}-cv\sim \tau (1-\gamma(x))v.
$$
Since all the points in $\Lambda_{\partial D}(p)$ attains the minimum of the function: $\partial D\ni x\longmapsto \vert x-p\vert$,
roughly speaking, one may expect, as $\tau\longrightarrow\infty$
$$\displaystyle
J(\tau)\sim\tau\int_{\partial D}(1-\gamma)v^2dS
$$
and
$$\displaystyle
J(\tau)+\int_{\partial D}
\frac{1}{c}\left\vert\frac{\partial v}{\partial\nu}-cv\right\vert^2dS
=\int_{\partial D}\frac{1}{c}\left(\frac{\partial v}{\partial\nu}-cv\right)\frac{\partial v}{\partial\nu}dS
\sim \tau\int_{\partial D}\frac{1-\gamma}{\gamma}v^2dS.
$$
These suggest (i) and (ii) of Lemma 2.3.  Note also that
$1-\gamma(x)\le (1-\gamma(x))/\gamma(x)$ if $\gamma(x)>0$.

Now it is easy to see that from Lemmas 2.2 and 2.3 one obtains ($\ref{eq:2.5}$)
and other statements of Theorem 2.1.
In the proof we never make use of the idea of {\it geometrical optics} which is classical.
Everything can be done in the context of the weak solution of \cite{DL} and main tool is just integration by parts.
Note that in Theorem 2.1 $\gamma$ is just essentially bounded on $\partial D$ and thus may have, for example,
a first kind of discontinuity.

\subsection{Curvatures and counting number}

Let $z\in\Bbb R^3$ and $0<r$.  In what follows we denote by $B_r(z)$ the open ball centred at $z$ and
with radius $r$.

Let $q\in\partial D$.  Given $v\in T_q(\partial D)$
define $\displaystyle
S_q(\partial D)v=-\frac{d}{dt}(\nu_{q(t)})\vert_{t=0}$,
where $q(t)\in\partial D$, $q(0)=q$ and $dq/dt(0)=v$.
We have $S_q(\partial D)v\in T_q(\partial D)$.
The operator $S_q(\partial D): T_q(\partial D)\longrightarrow T_q(\partial D)$ is called the {\it shape
operator} (or {\it Weingarten map}) of $\partial D$ at $q$ derived from $\nu$.
The shape operator is symmetric with respect to
the induced inner product on $T_q(\partial D)$ and its eigenvalues $k_1(q)\le k_2(q)$
are called the {\it principle curvatures} at $q$.
$K_{\partial D}(q)=k_1(q)k_2(q)$ and
$H_{\partial D}(q)=(k_1(q)+k_2(q))/2$ are called the Gauss and mean curvatures at $q$,
respectively.

Let $p\in\Bbb R^3\setminus\overline D$.
Let $q'\in \partial B_{d_{\partial D}(p)}(p)$ and $S_{q'}(\partial B_{d_{\partial D}(p)}(p))$ denote the shape operator of
$\partial B_{d_{\partial D}(p)}(p)$ at $q'$ derived from the unit {\it inward} normal to $\partial B_{d_{\partial D}(p)}(p)$.
If $q\in\Lambda_{\partial D}(p)$, then we have $q\in \partial B_{d_{\partial D}(p)}(p)$,
$T_q(\partial D)=T_q(\partial B_{d_{\partial D}(p)}(p))$
and $S_q(\partial B_{d_{\partial D}(p)}(p))-S_q(\partial D)\ge 0$.
Since $\displaystyle
S_q(\partial B_{d_{\partial D}(p)}(p))=(1/d_{\partial D}(p))I$,
we have
%$$\displaystyle
%\text{det}\,(S_q(\partial B_{d_{\partial D}(p)}(p))-S_q(\partial D))
%=(\lambda-k_1(q))(\lambda-k_2(q)),
%\tag {2.20}
%$$
\begin{eqnarray}\label{eq:2.20}
\displaystyle
\text{det}\,(S_q(\partial B_{d_{\partial D}(p)}(p))-S_q(\partial D))
=(\lambda-k_1(q))(\lambda-k_2(q)),
\end{eqnarray}
where $\lambda=1/d_{\partial D}(p)$.

\proclaim{\noindent Theorem 2.2(\cite{IWCOE})).}
Let $\gamma\equiv 0$.
Assume that $\partial D$ is $C^3$ and $\beta\in C^2(\partial D)$; $\Lambda_{\partial D}(p)$ is finite and satisfies
%$$\displaystyle
%\text{det}\,(S_q(\partial B_{d_{\partial D}(p)}(p))-S_q(\partial D))>0,\,\,\forall q\in\Lambda_{\partial D}(p).
%\tag {2.21}
%$$
\begin{eqnarray}\label{eq:2.21}
\displaystyle
\text{det}\,(S_q(\partial B_{d_{\partial D}(p)}(p))-S_q(\partial D))>0,\,\,\forall q\in\Lambda_{\partial D}(p).
\end{eqnarray}
If $T>2\text{dist}\,(D,B)$,
then we have
%$$\begin{array}{c}
%\displaystyle
%\lim_{\tau\longrightarrow\infty}
%\tau^4e^{2\tau\text{dist}\,(D,B)}
%I_B(\tau)
%=\frac{\pi}{2}
%\left(\frac{\eta}{d_{\partial D}(p)}\right)^2\mbox{\boldmath $A$}_{\partial D}(p),
%\end{array}
%\tag {2.22}
%$$
\begin{eqnarray}\label{eq:2.22}
\displaystyle
\lim_{\tau\longrightarrow\infty}
\tau^4e^{2\tau\text{dist}\,(D,B)}
I_B(\tau)
=\frac{\pi}{2}
\left(\frac{\eta}{d_{\partial D}(p)}\right)^2\mbox{\boldmath $A$}_{\partial D}(p),
\end{eqnarray}
where
$$\begin{array}{c}
\displaystyle
\mbox{\boldmath $A$}_{\partial D}(p)
=\sum_{q\in\Lambda_{\partial D}(p)}
\frac{1}{\displaystyle
\sqrt{
\text{det}\,(S_q(\partial B_{d_{\partial D}(p)}(p))-S_q(\partial D))}}
\end{array}
$$

\endproclaim

Using Theorem 2.2, one can give a procedure for extracting the curvatures at a known first refelection point.
More precisely, let $p\in\Bbb R^3\setminus\overline D$ and $q\in\Lambda_{\partial D}(p)$.
From ($\ref{eq:2.20}$) we have
%$$\displaystyle
%\text{det}\,(S_q(\partial B_{d_{\partial D}(p)}(p))-S_q(\partial D))=Q(\lambda)
%\equiv
%\lambda^2-2H_{\partial D}(q)\lambda+K_{\partial D}(q),
%\tag {2.23}
%$$
\begin{eqnarray}\label{eq:2.23}
\displaystyle
\text{det}\,(S_q(\partial B_{d_{\partial D}(p)}(p))-S_q(\partial D))=Q(\lambda)
\equiv
\lambda^2-2H_{\partial D}(q)\lambda+K_{\partial D}(q),
\end{eqnarray}
where $\lambda=1/d_{\partial D}(p)$.
Replace $\displaystyle
p\longrightarrow p_j=p-s_j\nu_q,\,\, j=1,2,\,\,0<s_1<s_2<d_{\partial D}(p)$.
Then $\Lambda_{\partial D}(p_j)=\{q\}$ and
$\text{det}\,(S_q(\partial B_{d_{\partial D}(p_j)}(p_j))-S_q(\partial D))>0$
since $\displaystyle
S_q(\partial B_{d_{\partial D}(p_j)}(p_j))>S_q(\partial B_{d_{\partial D}(p)}(p))$
and $\displaystyle
S_q(\partial B_{d_{\partial D}(p)}(p))-S_q(\partial D)\ge 0\,\,(\text{$q$ attains $\min_{x\in\partial D}\vert x-p\vert$})$.

Let $B_1$ and $B_2$ denote two open balls cetred at $p-s_j\nu_q$, $j=1,2$, respectively with
$0<s_1<s_2<d_{\partial D}(p)$ and satisfy $\overline B_1\cup\overline B_2\subset\Bbb R^3\setminus\overline D$.
Let $T>2\max\,\text{dist}\,(D,B_j)$ and $f=\chi_{B_j}$.
Applying ($\ref{eq:2.22}$) to this case, we obtain
$$\displaystyle
\lim_{\tau\longrightarrow\infty}
\tau^4e^{2\tau\text{dist}\,(D,B_j)}
I_{B_j}(\tau)
=\frac{\pi}{2}
\left(\frac{\text{diam}\,B_j}{2d_{\partial D}(p_j)}\right)^2
\frac{1}{\sqrt{Q(\lambda_j)}},
$$
where $\lambda_j=1/d_{\partial D}(p_j)$.
Since $\text{dist}\,(D,B_j)=d_{\partial D}(p_j)-s_j$ and $d_{\partial D}(p_j)=d_{\partial D}(p)-s_j$,
one can know $Q(\lambda_j)$ with $j=1,2$ from $u_f(x,t)$ given at all $(x,t)\in B_j\times\,]0,\,T[$
for $f=\chi_{B_j}$ with $j=1,2$.

Then, solving the system
%$$\displaystyle
%\left(\begin{array}{cc}
%\displaystyle -2\lambda_1 & 1\\
%\\
%\displaystyle
%-2\lambda_2 & 1
%\end{array}
%\right)
%\left(\begin{array}{c}
%\displaystyle
%H_{\partial D}(q)\\
%\\
%\displaystyle
%K_{\partial D}(q)
%\end{array}
%\right)
%=\left(\begin{array}{c}
%\displaystyle
%Q(\lambda_1)-\lambda_1^2\\
%\\
%\displaystyle
%Q(\lambda_2)-\lambda_2^2
%\end{array}
%\right),
%\tag {2.24}
%$$
\begin{eqnarray}\label{eq:2.24}
\displaystyle
\left(\begin{array}{cc}
\displaystyle -2\lambda_1 & 1\\
\\
\displaystyle
-2\lambda_2 & 1
\end{array}
\right)
\left(\begin{array}{c}
\displaystyle
H_{\partial D}(q)\\
\\
\displaystyle
K_{\partial D}(q)
\end{array}
\right)
=\left(\begin{array}{c}
\displaystyle
Q(\lambda_1)-\lambda_1^2\\
\\
\displaystyle
Q(\lambda_2)-\lambda_2^2
\end{array}
\right),
\end{eqnarray}
we obtain both $K_{\partial D}(q)$ and $H_{\partial D}(q)$.  Thus, one can know an approximate shape of the obstacle
in a neighbourhood of $q$.  Note that if $d_{\partial D}(p)\longrightarrow\infty$, then $\lambda_j\longrightarrow 0$
and thus it will be difficult to extract $H_{\partial D}(q)$ from ($\ref{eq:2.24}$).

Another simple corollay is a formula for counting the number of unknown spherical obstacles with
the same and known radius nearest to the center of the support of $f$.
Assume that $\displaystyle
D=B_{\epsilon}(x_1)\cup\cdots\cup B_{\epsilon}(x_m)$,
where $B_{\epsilon}(x_j)$, $j=1,\cdots,m$ is the open ball centred at $x_j$ with a known radius $\epsilon>0$
and $\overline{B_{\epsilon}(x_i)}\cap\overline{B_{\epsilon}(x_j)}=\emptyset$ if $i\not=j$.

Given $p\in\Bbb R^3\setminus\overline D$ it is easy to see that:
$\Lambda_{\partial D}(p)$ consists of finite points; ($\ref{eq:2.21}$) is satisfied;
there exists at most one first reflection point going from $p$ on each $\partial B_{\epsilon}(x_j)$.
Therefore, one can apply Theorem 2.2 to this case and obtain the formula which
enables us to know the counting number of the balls which are closest to the centre of $B$, that is,
$$\begin{array}{c}
\displaystyle
\sharp\Lambda_{\partial D}(p)
=\left(\frac{1}{d_{\partial D}(p)}+\frac{1}{\epsilon}\right)
\frac{2}{\pi}\left(\frac{\text{diam}\,B_{d_{\partial D}(p)}(p)}{\text{diam}\,B}\right)^2
\lim_{\tau\longrightarrow\infty}\tau^4e^{2\tau\text{dist}\,(D,B)}
I_B(\tau),
\end{array}
$$
where
$$\displaystyle
\Lambda_{\partial D}(p)
=\left\{x_i+\epsilon\frac{p-x_i}{\vert p-x_i\vert}\,\vert\,\vert p-x_i\vert=\min_{j}\vert p-x_j\vert\right\}.
$$

\subsection{A sketch of the proof of Theorem 2.2}
Let $\gamma\equiv 0$.
Integration by parts yields
$$\displaystyle
J(\tau)=\int_{D}(\vert\nabla v\vert^2+\tau^2\vert v\vert^2)dx-\int_{\partial D}\beta \vert v\vert^2dS.
$$
Applying a trace theorem \cite{G} to the second integral on this right-hand side, we see tha
$J(\tau)>0$ for all $\tau>>1$.
Then, we have, as $\tau\longrightarrow\infty$
%$$\displaystyle
%E(\tau)=J(\tau)(1+O(\tau^{-1/2}))
%\tag {2.25}
%$$
\begin{eqnarray}\label{eq:2.25}
\displaystyle
E(\tau)=J(\tau)(1+O(\tau^{-1/2}))
\end{eqnarray}
and thus from ($\ref{eq:2.13}$) we obtain
%$$\displaystyle
%I_B(\tau)=2J(\tau)(1+O(\tau^{-1/2}))+O(\tau^{-1}e^{-\tau T}).
%\tag {2.26}
%$$
\begin{eqnarray}\label{eq:2.26}
\displaystyle
I_B(\tau)=2J(\tau)(1+O(\tau^{-1/2}))+O(\tau^{-1}e^{-\tau T}).
\end{eqnarray}
Using the Laplace method \cite{BH}, one can expand $J(\tau)$  under the condition ($\ref{eq:2.21}$) and we find its leading term
which contains information about the geometry of $\partial D$ at all the first reflection points, going from the centre
of $B$ to $\partial D$.   This yields ($\ref{eq:2.22}$).

Thus the crucial point of the proof of Theorem 2.2 is the derivation of ($\ref{eq:2.25}$).
It is a combination of a modification of the Lax-Phillips reflection argument in \cite{LP}
and a change of a dependent variable near $\partial D$.
Here we describe the idea of the derivation of ($\ref{eq:2.25}$) in the simplest case $\gamma=\beta\equiv 0$.

Since $G\equiv 0$, it follows from ($\ref{eq:2.8}$)
$$\displaystyle
E(\tau)
=\int_{\partial D}\frac{\partial v}{\partial\nu}RdS-e^{-\tau T}\int_{\Bbb R3\setminus\overline D}FRdx
$$
and applying the boundary condition in ($\ref{eq:2.6}$) to $J(\tau)$, we obtain
$$\displaystyle
E(\tau)-J(\tau)
=\int_{\partial D}\left(\frac{\partial v}{\partial\nu}R+\frac{\partial R}{\partial\nu}v\right)dS-e^{-\tau T}\int_{\Bbb R3\setminus\overline D}FRdx.
$$
Choose $\tilde{v}(x), x\in\Bbb R^3\setminus D$ in such a way that $\tilde{v}$ has a compact support
and satisfies
%$$\left\{
%\begin{array}{l}
%\displaystyle
%\tilde{v}=v\,\,\text{on}\,\partial D,\\
%\\
%\displaystyle
%\frac{\partial\tilde{v}}{\partial\nu}=-\frac{\partial v}{\partial\nu}\,\,\text{on}\,\partial D.
%\end{array}
%\right.
%\tag {2.27}
%$$
\begin{eqnarray}\label{eq:2.27}
\left\{
\begin{array}{ll}
\displaystyle
\tilde{v}=v & \text{on}\quad\partial D,
\\[10pt]
\displaystyle
\frac{\partial\tilde{v}}{\partial\nu}=-\frac{\partial v}{\partial\nu} & \text{on}\quad\partial D.
\end{array}
\right.
\end{eqnarray}
Then, integration by parts and ($\ref{eq:2.6}$) gives
$$\displaystyle
\int_{\partial D}\left(\frac{\partial v}{\partial\nu}R+\frac{\partial R}{\partial\nu}v\right)dS
=\int_{\Bbb R^3\setminus\overline D}(\triangle-\tau^2)\tilde{v}\cdot Rdx
-e^{-\tau T}\int_{\Bbb R^3\setminus\overline D}\tilde{v}Fdx.
$$
Hence we obtain
%$$\begin{array}{c}
%\displaystyle
%E(\tau)-J(\tau)
%=\int_{\Bbb R^3\setminus\overline D}(\triangle-\tau^2)\tilde{v}\cdot Rdx
%-e^{-\tau T}\int_{\Bbb R^3\setminus\overline D}F(R+\tilde{v})dx.
%\end{array}
%\tag {2.28}
%$$
\begin{eqnarray}\label{eq:2.28}
\displaystyle
E(\tau)-J(\tau)
=\int_{\Bbb R^3\setminus\overline D}(\triangle-\tau^2)\tilde{v}\cdot Rdx
-e^{-\tau T}\int_{\Bbb R^3\setminus\overline D}F(R+\tilde{v})dx.
\end{eqnarray}
The point is the choice of $\tilde{v}$.
Let $x^r$ denote the {\it reflection} in the tubular neighbourhood $\{x\in\Bbb R^3\setminus D\,\vert d_{\partial D}(x)<2\delta_0\}$
of $\partial D$ with sufficiently small $\delta_0>0$.  It is given by $x^r=2q(x)-x$, where
$q(x)$ denote the unique point on $\partial D$ such that $d_{\partial D}(x)=\vert x-q(x)\vert$.
It is known that $q(x)$ is $C^2$ for $x\in\Bbb R^3\setminus D$ with $d_{\partial D}(x)<2\delta_0$ if $\partial D$
is $C^3$ (\cite{GT}).
Choose a cutoff function $\phi_{\delta}$ with $0<\delta<\delta_0$ which satisfies $0\le\phi_{\delta}(x)\le 1$;
$\phi_{\delta}(x)=1$ if $d_{\partial D}(x)<\delta$; $\phi_{\delta}(x)=0$ if $d_{\partial D}(x)>2\delta$; $\vert\nabla\phi_{\delta}(x)\vert
\le C\delta^{-1}$; $\vert\nabla^2\phi_{\delta}(x)\vert\le C\delta^{-2}$.

Define
$$
\displaystyle\tilde{v}(x)=\phi_{\delta}(x)v(x^r).
$$
Clearly ($\ref{eq:2.27}$) is satisfied with this $\tilde{v}$.
A direct computation gives
%$$\displaystyle
%(\triangle-\tau^2)\tilde{v}(x)
%=\phi(x)d_{\partial D}(x)\sum_{i,j}a_{ij}(x)(\partial_i\partial_jv)(x^r)
%+(\text{lower order terms}),
%\tag {2.29}
%$$
\begin{eqnarray}\label{eq:2.29}
\displaystyle
(\triangle-\tau^2)\tilde{v}(x)
=\phi_{\delta}(x)d_{\partial D}(x)\sum_{i,j}a_{ij}(x)(\partial_i\partial_jv)(x^r)
+(\text{lower order terms}),
\end{eqnarray}
where $a_{ij}(x)$ with $i,j=1,2,3$ are $C^1$ for $x\in\Bbb R^3\setminus D$
with $d_{\partial D}(x)<2\delta_0$ and independent of $\tau$, $\phi_{\delta}$ and $v$.
Note that the computation is based on the formula
$$\displaystyle
2q'(x)-I
=I-2\nu_{q(x)}\otimes\nu_{q(x)}
-2d_{\partial D}(x)(\nu_{q(x)})',
$$
where $x\in\Bbb R^3\setminus\overline D$ and $d_{\partial D}(x)<<1$;
$q'(x)$ denotes the Jacobian matrix of the map: $x\longmapsto q(x)$.
It is a consequence of the expression
$q(x)=x-d_{\partial D}(x)\nu_{q(x)}$ and the formula $\nabla(d_{\partial D}(x))=\nu_{q(x)}$.

The point is $d_{\partial D}(x)$ in the first term on the right-hand side of ($\ref{eq:2.29}$).
By using the change of variable $x=y^r$ we have
%$$\begin{array}{c}
%\displaystyle
%\int_{\Bbb R^3\setminus\overline D}(\triangle-\tau^2)\tilde{v}\cdot Rdx\\
%\\
%\displaystyle
%=\sum_{i,j}\int_{D}\phi(y^r)d_{\partial D}(y^r)a_{ij}(y^r)(\partial_i\partial_j v)(y)R(y^r)J(y)dy+\cdots,
%\end{array}
%\tag {2.30}
%$$
\begin{eqnarray}\label{eq:2.30}
\begin{array}{c}
\displaystyle
\int_{\Bbb R^3\setminus\overline D}(\triangle-\tau^2)\tilde{v}\cdot Rdx
\\[10pt]
\displaystyle
=\sum_{i,j}\int_{D}\phi_{\delta}(y^r)d_{\partial D}(y^r)a_{ij}(y^r)(\partial_i\partial_j v)(y)R(y^r)J(y)dy+\cdots,
\end{array}
\end{eqnarray}
where $J(y)$ denote the Jacobian of the map: $y\longmapsto y^r$.  Since $d_{\partial D}(y^r)\equiv d_{\partial D}(y)=0$
on $\partial D$, integration by parts yields
$$\begin{array}{c}
\displaystyle
\int_{D}\phi_{\delta}(y^r)d_{\partial D}(y^r)a_{ij}(y^r)(\partial_i\partial_j v)(y)R(y^r)J(y)dy\\
\\
\displaystyle
=-\int_{D}\partial_{i}\{\phi_{\delta}(y^r)d_{\partial D}(y)a_{ij}(y^r)R(y^r)J(y)\}\partial_jv(y)dy.
\end{array}
$$
Hereafter simply estimating this right-hand side together with other terms in ($\ref{eq:2.30}$), we obtain
$$\begin{array}{c}
\displaystyle
\left\vert\int_{\Bbb R^3\setminus\overline D}(\triangle-\tau^2)\tilde{v}\cdot Rdx\right\vert\\
\\
\displaystyle
\le C((\delta\Vert\nabla R^r\Vert_{L^2(D_{\delta})}+\delta^{-1}\Vert R^r\Vert_{L^2(D_{\delta})})\Vert\nabla v\Vert_{L^2(D)}
+\delta^{-2}\Vert R^r\Vert_{L^2(D_{\delta})}\Vert v\Vert_{L^2(D)}),
\end{array}
$$
where $D_{\delta}=\{y\in D\,\vert\,d_{\partial D}(y)<2\delta\}$ and $R^r(y)=R(y^r)$.

Here we note that
$$\displaystyle
\Vert\nabla v\Vert_{L^2(D)}\le J(\tau)^{1/2},\,\,
\Vert v\Vert_{L^2(D)}\le\tau^{-1}J(\tau)^{1/2}
$$
and
$$
\displaystyle
\Vert\nabla R^r\Vert_{L^2(D_{\delta})}
\le CE(\tau)^{1/2},\,\,
\Vert R^r\Vert_{L^2(D_{\delta})}\le
\tau^{-1}E(\tau)^{1/2}.
$$
Choosing $\delta=\tau^{-1/2}$ with $\tau>>1$,
we finally obtain
$$\displaystyle
\left\vert\int_{\Bbb R^3\setminus\overline D}(\triangle-\tau^2)\tilde{v}\cdot Rdx\right\vert
\le C\tau^{-1/2}(E(\tau)J(\tau))^{1/2}.
$$
From these together with ($\ref{eq:2.28}$) and the estimate
$\displaystyle
\Vert \tilde{v}\Vert_{L^2(\Bbb R^3\setminus\overline D)}\le C\Vert v\Vert_{L^2(D)}$,
we obtain
%$$\displaystyle
%\vert E(\tau)-J(\tau)\vert\le C(\tau^{-1/2}(E(\tau)J(\tau))^{1/2}+e^{-\tau T}E(\tau)^{1/2}+e^{-\tau T}J(\tau)^{1/2})
%\tag {2.31}
%$$
\begin{eqnarray}\label{eq:2.31}
\displaystyle
\vert E(\tau)-J(\tau)\vert\le C(\tau^{-1/2}(E(\tau)J(\tau))^{1/2}+e^{-\tau T}E(\tau)^{1/2}+e^{-\tau T}J(\tau)^{1/2})
\end{eqnarray}
and hence
$$
\displaystyle
(1-2C\tau^{-1}-2e^{-\tau T})E(\tau)
\le(1+2C+2e^{-\tau T})J(\tau)+4e^{-\tau T}.
$$
Therefore, there exist potive constants $C'$ and $\tau_0$ such that, for all $\tau\ge\tau_0$
$$\displaystyle
E(\tau)\le C'(J(\tau)+e^{-\tau T}).
$$
From the case (i) of Lemma 2.3 one can conclude that
$e^{-\tau T}/J(\tau)$ is decreasing as $\tau\longrightarrow\infty$ provided
$T>2\text{dist}\,(D,B)$.  Thus we have $E(\tau)=O(J(\tau))$ as $\tau\longrightarrow\infty$.
Applying this together with the trivial estimate $J(\tau)=O(1)$ as $\tau\longrightarrow\infty$ to the
right-hand side on ($\ref{eq:2.31}$), we finally obtain ($\ref{eq:2.25}$).

\subsection{Quantitative state of the surface}

A combination of ($\ref{eq:2.26}$) and the second order term of the asymptotic expansion of $J(\tau)$ as $\tau\longrightarrow\infty$ yields
information about the value of $\beta$ at all the first reflection points, going from the centre to $\partial D$.

\proclaim{\noindent Theorem 2.3(\cite{IWCOE}).}  Let $\gamma\equiv 0$.
Assume that $\partial D$ is $C^5$ and $\beta\in C^2(\partial D)$;
$\Lambda_{\partial D}(p)$ is finite and satisfies ($\ref{eq:2.21}$).
For each $q\in \Lambda_{\partial D}(p)$ let $\mbox{\boldmath
$e$}_j$, $j=1,2$ be an orthonormal basis of the tangent space at $q$
of $\partial D$ with $\mbox{\boldmath $e$}_1\times\mbox{\boldmath
$e$}_2=\nu_q$. Choose an open ball $U$ centred at $q$ with radius
$r_q$ in such a way that there exist a $h\in C^5_0(\Bbb R^2)$
with $h(0,0)=0$ and $\nabla h(0,0)=0$ such that $U\cap\partial
D=\{q+\sigma_1\mbox{\boldmath $e$}_1+\sigma_2\mbox{\boldmath
$e$}_2+h(\sigma_1,\sigma_2)\nu_q\,\vert\,
\sigma_1^2+\sigma_2^2+h(\sigma_1,\sigma_2)^2<r_q^2\}$.

If $T>2\text{dist}\,(D,B)$,
then we have
%$$\begin{array}{c}
%\displaystyle
%\lim_{\tau\longrightarrow\infty}
%\tau^5\left\{e^{2\tau\text{dist}\,(D,B)}
%I_B(\tau)-\frac{1}{\tau^4}\frac{\pi}{2}
%\left(\frac{\eta}{d_{\partial D}(p)}\right)^2\mbox{\boldmath $A$}_{\partial D}(p)\right\}\\
%\\
%\displaystyle
%=-\frac{\pi\eta}{d_{\partial D}(p)^2}\mbox{\boldmath $A$}_{\partial D}(p)
%+\frac{\pi}{2}\eta^2\mbox{\boldmath $B$}_{\partial D}(p),
%\end{array}
%\tag {2.32}
%$$
\begin{eqnarray}\label{eq:2.32}
\begin{array}{c}
\displaystyle
\lim_{\tau\longrightarrow\infty}
\tau^5\left\{e^{2\tau\text{dist}\,(D,B)}
I_B(\tau)-\frac{1}{\tau^4}\frac{\pi}{2}
\left(\frac{\eta}{d_{\partial D}(p)}\right)^2\mbox{\boldmath $A$}_{\partial D}(p)\right\}
\\[12pt]
\displaystyle
=-\frac{\pi\eta}{d_{\partial D}(p)^2}\mbox{\boldmath $A$}_{\partial D}(p)
+\frac{\pi}{2}\eta^2\mbox{\boldmath $B$}_{\partial D}(p),
\end{array}
\end{eqnarray}
where
$$\begin{array}{c}
\displaystyle
\mbox{\boldmath $B$}_{\partial D}(p)
=\sum_{q\in\Lambda_{\partial D}(p)}
\frac{\mbox{\boldmath $C$}_{\partial D}(q)}
{\displaystyle
\sqrt{\text{det}\,(S_q(\partial B_{d_{\partial D}(p)}(p))-S_q(\partial D))}},
\end{array}
$$
$$\begin{array}{c}
\displaystyle
\mbox{\boldmath $C$}_{\partial D}(q)
=-\frac{1}{d_{\partial D}(p)^3}+\frac{11-12d_{\partial D}(p)H_{\partial D}(q)}
{\displaystyle
8d_{\partial D}(p)^5
\text{det}\,(S_q(\partial B_{d_{\partial D}(p)}(p))-S_q(\partial D))}\\
\\
\displaystyle
-\frac{1}{4d_{\partial D}(p)^2}
h_{\sigma_p\sigma_q\sigma_r}(0)h_{\sigma_s\sigma_t\sigma_u}(0)
\left(\frac{1}{4}B_{ps}B_{qr}B_{tu}+\frac{1}{6}B_{ps}B_{qt}B_{ru}\right)\\
\\
\displaystyle
+\frac{1}{16d_{\partial D}(p)^2}
h_{\sigma_p\sigma_q\sigma_r\sigma_s}(0)B_{pr}B_{qs}
-\frac{\beta(q)}{d_{\partial D}(p)^2}
\end{array}
$$
and
$$\displaystyle
(B_{pq})
=-\left(\frac{1}{d_{\partial D}(p)}I_2-\nabla^2h(0)\right)^{-1}.
$$

\endproclaim

Note that we have used the summation convention where repeated indicies are to be summed from $1$ to $2$.
The explicit second-order term of an expansion of Laplace type integral $J(\tau)$ is essential and it
is an application of an expansion formula of the Laplace type integral.  For this see \cite{BH}.

As a corollary of Theorem 2.3 we obtain a procedure for calculating the value of $\beta$ at a known point on $\Lambda_{\partial D}(p)$
provided $\gamma\equiv 0$.
More precisely, we assume that:

\noindent
(i)  we know in advance a point $q\in \Lambda_{\partial D}(p)$;

\noindent
(ii)  we know that $\partial D$ near $q$
is given by making a {\it rotation} around the normal at $q$ of a
{\it graph} of a function $h$ defined on the tangent plane at $q$ of $\partial D$ and
that, in an appropriate orthogonal coordinates on the tangent
plane, say $\sigma=(\sigma_1,\sigma_2)$, the Taylor expansion of
the function at $\sigma=0$ has the form
$h(\sigma_1,\sigma_2)=\sum_{2\le \vert\alpha\vert\le
4}h_{\alpha}\sigma^{\alpha}+\cdots$ with {\it known} coefficients
$h_{\alpha}$ for $2\le\vert\alpha\vert\le 4$.

Note that, from (i) we know also $d_{\partial D}(p)=\vert p-q\vert$, $\nu_q=(p-q)/\vert p-q\vert$ and
the tangent plane $(x-q)\cdot\nu_q=0$ at $q$ of $\partial D$.

Fix $s\in]0,\,d_{\partial D}(p)[$.
Choose an open ball $B'$ centred at $p-s\nu_q$ and satisfying $\overline{B'}\subset B_{d_{\partial D}(p)}(p)$.
Let $T>2\text{dist}\,(D,\,B')$.
Generate the wave $u_f$ by $f=\chi_{B'}$ and observe the wave on $B'$ over time interval $]0,\,T[$.
Since $\Lambda_{\partial D}(p-s\nu_q)=\{q\}$ and ($\ref{eq:2.21}$) for $p$ replaced with $p-s\nu_q$ is satisfied,
one gets ($\ref{eq:2.32}$) in which ball $B$ is replaced with $B'$ and $p$ replaced with $p-s\nu_q$.
Therefore, we obtain $\mbox{\boldmath $C$}_{\partial D}(q)$ which yields a linear equation with unknown
$\beta(q)$ and thus solving this, one obtains $\beta(q)$.  Note that, in this procedure we do not assume
that $\Lambda_{\partial D}(p)$ is finite.

Now it is natural to consider the following problem.

{\bf\noindent Open problem 2.1.}  Asuume that, say, $\gamma$ is sufficiently smooth on $\partial D$.
Find a formula for calculating $\gamma$ at a known $q\in\Lambda_{\partial D}(p)$ from $u_f$ on $B'\times\,]0,\,T[$
generated by $f=\chi_{B'}$, where $B'$ is the same as above.

\noindent
The point is to find the asymptotic profile of $E(\tau)$ in ($\ref{eq:2.13}$) as $\tau\longrightarrow\infty$ in terms of $v$ on $D$.
For one-space dimensional case we have an explicit formula.  See \cite{IW2}.

\subsection{Other wave equations}

In \cite{IW2} a result analogous to Theorem 2.1 has been established also for the equation
$$
\displaystyle
\alpha(x)\partial_t^2 u-\triangle u=0\,\,\text{in}\,\Bbb R^3\times\,]0,\,T[
$$
provided: $\eta\ge \alpha(x)\ge\eta^{-1}$ a.e. $x\in\Bbb R^3$ for a positive constant $\eta$ and
$\alpha(x)=1$ a.e. $x\in\Bbb R^3\setminus D$; there exists a positive constant $C$ such that $\alpha(x)\le 1-C$
a.e. $x\in D$ or $\alpha(x)\ge 1+C$ a.e. $x\in D$.

In \cite{IW} the equation
$$\displaystyle
\partial_t^2u-\nabla\cdot A(x)\nabla u=0\,\,\text{in}\,\Bbb R^3\times\,]0,\,T[
$$
with a $3\times 3$ uniformly positive definite real symmetric matrix-valued function coefficient $A(x)$ satisfying $A(x)=I_3$ a.e. $x\in\Bbb R^3\setminus D$
has been studied.  It is assumed that each component of $A(x)$ is essentially bounded and
there exists a positive constant $C$
such that $(A(x)-I_3)\xi\cdot\xi\ge C\vert\xi\vert^2$ a.e. $x\in D$ and all
$\xi\in\Bbb R^3$ or $-(A(x)-I_3)\xi\cdot\xi\ge C\vert\xi\vert^2$ a.e. $x\in D$
and all $\xi\in\Bbb R^3$.  Then, it is clear that Theorem 1.2 in \cite{IW} for this equation yields also a result analogous to Theorem 2.1.
However, for both equations there is no result corresponding to Theorems 2.2 and 2.3 via the Enclosure Method.

\subsection{Interior problem in time domain}

Let $D$ be a bounded domain of $\Bbb R^3$ with $C^2$-boundary.
Given $f\in L^2(D)$ satisfying $\text{supp}\,f\subset D$
denote by $u_f$ the solution of the following initial boundary value problem for the wave equation:
%$$\left\{
%\begin{array}{c}
%\displaystyle
%\partial_t^2u-\triangle u=0\,\,\text{in}\, D\times\,]0,\,T[,\\
%\\
%\displaystyle
%-\frac{\partial u}{\partial\nu}-\gamma(x)\partial_tu-\beta(x) u=0\,\,\text{on}\,\partial D\times\,]0,\,T[,
%\\
%\\
%\displaystyle
%u(x,0)=0\,\,\text{in}\,D,\,\,
%\partial_tu(x,0)=f(x)\,\,\text{in}\,D,
%\end{array}
%\right.
%\tag {2.33}
%$$
\begin{eqnarray}\label{eq:2.33}
\left\{
\begin{array}{ll}
\displaystyle
\partial_t^2u-\triangle u=0 & \text{in}\quad D\times\,]0,\,T[,
\\[10pt]
\displaystyle
-\frac{\partial u}{\partial\nu}-\gamma(x)\partial_tu-\beta(x) u=0 & \text{on}\quad\partial D\times\,]0,\,T[,
\\[10pt]
\displaystyle
u(x,0)=0 & \text{in}\quad D,
\\[10pt]
\displaystyle
\partial_tu(x,0)=f(x) & \text{in}\quad D,
\end{array}
\right.
\end{eqnarray}
where $\nu$ denotes the unit {\it outward normal} to $D$ on $\partial D$,
$\beta$ and $\gamma$ are the same as those in ($\ref{eq:2.1}$).

\noindent
{\bf Problem 2.2.}  Let $B$ be an open ball and satisfy $\overline B\subset D$.
Generate $u=u_f$ by the initial data $f=\chi_B$ and observe $u$ on $B$ over
time interval $]0,\,T[$.  Extract information about the geometry of $\partial D$, $\gamma$ and $\beta$ from the observed data.

In \cite{IWcavity}, we considered the case when $\gamma=\beta=0$ in ($\ref{eq:2.33}$) and obtained two theorems corresponding to Theorems 2.1
and 2.2. It will be possible to obtain a theorem corresponding to Theorem 2.1 for general $\gamma$ and $\beta$;
theorems corresponding to Theorems 2.2 and 2.3 for $\gamma=0$ and general $\beta$.
Thus, a real problem to be solved should be the same as Open problem 2.1.

\section{Further applications and problems}

\subsection{Bistatic data, spheroid and simultaneous rotation}

Let $0<T<\infty$.
Let $f\in L^2(\Bbb R^3)$ satisfy $\text{supp}\,f\cap\overline D=\emptyset$.
Let $u=u_f(x,t)$ be the solution of the initial boundary value problem:
%$$\left\{
%\begin{array}{c}
%\displaystyle
%\partial_t^2u-\triangle u=0\,\,\text{in}\,(\Bbb R^3\setminus\overline D)\times\,]0,\,T[,
%\\
%\\
%\displaystyle
%u=0\,\,\text{on}\,\partial D\times\,]0,\,T[,
%\\
%\\
%\displaystyle
%u(x,0)=0\,\,\text{in}\,\Bbb R^3\setminus\overline D,\,\,
%\partial_tu(x,0)=f(x)\,\,\text{in}\,\Bbb R^3\setminus\overline D.
%\end{array}
%\right.
%$$
$$\left\{
\begin{array}{ll}
\displaystyle
\partial_t^2u-\triangle u=0 & \text{in}\quad (\Bbb R^3\setminus\overline D)\times\,]0,\,T[,
\\[10pt]
\displaystyle
u=0 & \text{on}\quad \partial D\times\,]0,\,T[,
\\[10pt]
\displaystyle
u(x,0)=0 & \text{in}\quad \Bbb R^3\setminus\overline D,
\\[10pt]
\displaystyle
\partial_tu(x,0)=f(x) & \text{in}\quad\Bbb R^3\setminus\overline D.
\end{array}
\right.
$$

{\bf\noindent Problem 3.1.}
Let $B$ and $B'$ be two {\it known} open balls centred at $p\in\Bbb R^3$ and $p'\in\Bbb R^3$ with radii $\eta$
and $\eta'$, respectively such that $\overline B\cap\overline D=\emptyset$
and $\overline B'\cap \overline D=\emptyset$.
Let $\chi_B$ denote the characteristic function of $B$ and set $f=\chi_B$.
Assume that $D$ is {\it unknown}.
Extract information about the location and shape of $D$ from the data $u_f(x,t)$ given at all
$x\in B'$ and $t\in\,]0,\,T[$.

Let $\chi_{B'}$ denote the characteristic function of $B'$
and set $g=\chi_{B'}$.
The results of this subsection are concerned with the asymptotic
behaviour of the bistatic {\it indicator function} $I_{B,B'}$ defined by
$$\displaystyle
I_{B,B'}(\tau)=\int_{\Bbb R^3\setminus\overline D}(fv_g-w_fg)dx,
$$
where $w_f$ is the same as ($\ref{eq:2.4}$) and $v_g$ is the solution of ($\ref{eq:2.2}$) in which $f$ is replaced with $g$.
Note that $I_{B,B'}(\tau)$ can be computed from $w_f$ on $B$ and thus from $u_f$ on
$B\times\,]0,\,T[$.

Define $\displaystyle
\phi(x;y,y')
=\vert y-x\vert+\vert x-y'\vert,\,\,(x,y,y')\in\Bbb R^3\times\Bbb R^3\times\Bbb R^3$.
This is the length of the broken path connecting $y$ to $x$ and $x$ to $y'$.
We denote the {\it convex hull} of the set $F\subset\Bbb R^3$ by $[F]$.

\proclaim{\noindent Theorem 3.1(\cite{IW3}).}
Let $[\overline B\cup\overline B']\cap\partial D=\emptyset$
and $T$ satisfy
%$$\displaystyle
%T>\min_{x\in\partial D,\,y\in\partial B,\, y'\in\partial B'}\phi(x;y,y').
%\tag {3.1}
%$$
\begin{eqnarray}\label{eq:3.1}
\displaystyle
T>\min_{x\in\partial D,\,y\in\partial B,\, y'\in\partial B'}\phi(x;y,y').
\end{eqnarray}
Then, there exists a $\tau_0>0$ such that
$\displaystyle I_{B,B'}(\tau)>0$ for all $\tau\ge\tau_0$
and the formula
%$$\displaystyle
%\lim_{\tau\longrightarrow\infty}
%\frac{1}{\tau}
%\log
%I_{B,B'}(\tau)
%=-\min_{x\in\partial D,\,y\in\partial B,\, y'\in\partial B'}\phi(x;y,y'),
%\tag {3.2}
%$$
\begin{eqnarray}\label{eq:3.2}
\displaystyle
\lim_{\tau\longrightarrow\infty}
\frac{1}{\tau}
\log
I_{B,B'}(\tau)
=-\min_{x\in\partial D,\,y\in\partial B,\, y'\in\partial B'}\phi(x;y,y'),
\end{eqnarray}
is valid.

\endproclaim

It is easy to see that
$$\displaystyle
\min_{x\in\partial D,\, y\in \partial B,\,y'\in \partial B'}\phi(x;y,y')
=\min_{x\in\partial D}\phi(x;p,p')-(\eta+\eta').
$$
Thus formula ($\ref{eq:3.2}$) enables us to extract $\min_{x\in\partial D}\phi(x;p,p')$
from $u_f(x,t)$ given at all $x\in B'$ and $t\in]0,\,T[$.

The quantity $\min_{x\in\partial D}\phi(x;p,p')$ coincides with the shortest length
of the broken paths connecting $p$ to a point $q$ on $\partial D$ and $q$ to $p'$,
that is, the {\it first reflection distance} between $p$ and $p'$ by $D$.
Thus, Theorem 3.1 yields a mathematical method for extracting the first reflection distance
from the waveform of the observed wave.

Given $c>\vert p-p'\vert$ define $\displaystyle
E_c(p,p')=\{x\in\Bbb R^3\,\vert\,\phi(x;p,p')=c\}$.
This is a {\it spheroid} with focal points $p$ and $p'$.
Given direction $\omega\in S^2$ at $p'$ let $\zeta(\omega;p,p')$ denote the unique point on 
$E_c(p,p')\cap\{p'+s\omega\,\vert\,s>0\}$.
$\zeta(\omega;p,p,p')$ has the expression
$\displaystyle
\zeta(\omega;p,p')=p'+s(\omega;p,p')\omega$
with a unique $s(\omega;p,p')>0$ and the map $S^2:\omega\longmapsto \zeta(\omega;p,p')\in E_c(p,p')$ is bijective.

Define $\displaystyle
\Lambda_{\partial D}(p,p')=\{q\in\partial D\,\vert\,\min_{x\in\partial D}\phi(x;p,p')=\phi(q;p,p')\}$.
One can write $\Lambda_{\partial D}(p,p')=\partial D\cap E_c(p,p')$, where $c=\min_{x\in\partial D}\phi(x;p,p')$.

Similary to Section 2.2, using Theorem 3.1, one can make a decision
whether given direction $\omega\in S^2$ at $p'$ $\zeta(\omega;p,p')$
which is a point on $E_c(p,p')$ belongs to $\partial D$ or not.
It is based on the following chracterization of $\Lambda_{\partial D}(p,p')$.

\proclaim{\noindent Lemma 3.1(Proposition 5.1 in \cite{IW3}).}
Fix $s\in ]0,\,\eta'[$.
We have:

(i)  if $\zeta(\omega;p,p')$ belongs to $\partial D$, then $\min_{x\in\partial D}\phi(x;p,p'+s\omega)=c-s$;

(ii) if $\zeta(\omega;p,p')$ does not belong to $\partial D$, then $\min_{x\in\partial D}\phi(x;p,p'+s\omega)>c-s$.

\endproclaim

Therefore, we obtain the following characterization of $\Lambda_{\partial D}(p,p')$:
$$\displaystyle
\Lambda_{\partial D}(p,p')=\{\zeta(\omega;p,p')\,\vert\,\min_{x\in\partial D}\phi(x;p,p'+s\omega)=c-s\}.
$$

The procedure for finding $\zeta(\omega;p,p')$ belonging to $\partial D$ from a single set of the bistatic data
is the following:

Fix a large $T$ and $s\in]0,\,\eta'[$.

\noindent
{\bf Step 1.}   Generate $u_f$ by the initial data $f=\chi_B$ and observe $u_f$ on $B'$ over time interval $]0,\,T[$.

\noindent
{\bf Step 2.}   Choose an open ball $B''\subset B'$ centred at $p'+s\omega$.

\noindent
{\bf Step 3.}   Determine $\min_{x\in\partial D}\phi(x;p,p'+s\omega)$ from the restriction of $u_f$ in the first step
onto $B''\times\,]0,\,T[$ via Theorem 3.1.

From the computed value $\min_{x\in\partial D}\phi(x;p,p'+s\omega)$ in the third step,
one has: if $\min_{x\in\partial D}\phi(x;p,p'+s\omega)=c-s$, then $\zeta(\omega;p,p')$ belongs to $\partial D$;
if not, then $\zeta(\omega;p,p')$ does not belong to $\partial D$.

Therefore, in principle, one can determine all the points
in $\Lambda_{\partial D}(p,p')$ from $u_f$ on $B'\times\,]0,\,T[$ for $f=\chi_B$.
This is an advantage of the bistatic data not being seen in the monostatic data.

The next theoretical result is concerned with obtaining information about shape of $D$.
In the following theorem, for simplicity of description, we assume that $D$ is {\it convex}.
In this case $\Lambda_{\partial D}(p,p')$ consists of a single point $q=q(p,p')$.

\proclaim{\noindent Theorem 3.2(\cite{IW3}).}
Assume that $\partial D$ is $C^3$.
Let $c=\min_{x\in\partial D}\phi(x;p,p')$.
Let $T$ satisfy ($\ref{eq:3.1}$).
Then, we have
%$$\displaystyle
%\text{det}\,(S_{q(p,p')}(E_c(p,p'))-S_{q(p,p')}(\partial D))>0
%\tag {3.3}
%$$
\begin{eqnarray}\label{eq:3.3}
\displaystyle
\text{det}\,(S_{q(p,p')}(E_c(p,p'))-S_{q(p,p')}(\partial D))>0
\end{eqnarray}
and the formula
%$$\begin{array}{c}
%\displaystyle
%\lim_{\tau\longrightarrow\infty}
%\tau^4e^{\tau\min_{x\in\partial D,\,y\in\partial B,\, y'\in\partial B'}\phi(x;y,y')}
%I_{B,B'}(\tau)\\
%\\
%\displaystyle
%=\frac{\pi}{2}
%\left(\frac{\text{diam}\,B}{2\vert q-p\vert}\right)\cdot
%\left(\frac{\text{diam}\,B'}{2\vert q-p'\vert}\right)
%\cdot
%\frac{1}
%{\sqrt{\text{det}\,(S_{q}(E_c(p,p'))-S_{q}(\partial D))}}\vert_{q=q(p,p')},
%\end{array}
%\tag {3.4}
%$$
\begin{eqnarray}\label{eq:3.4}
\begin{array}{c}
\displaystyle
\lim_{\tau\longrightarrow\infty}
\tau^4e^{\tau\min_{x\in\partial D,\,y\in\partial B,\, y'\in\partial B'}\phi(x;y,y')}
I_{B,B'}(\tau)\\
\\
\displaystyle
=\frac{\pi}{2}
\left(\frac{\text{diam}\,B}{2\vert q-p\vert}\right)\cdot
\left(\frac{\text{diam}\,B'}{2\vert q-p'\vert}\right)
\cdot
\frac{1}
{\sqrt{\text{det}\,(S_{q}(E_c(p,p'))-S_{q}(\partial D))}}\vert_{q=q(p,p')},
\end{array}
\end{eqnarray}
is valid.

\endproclaim

\noindent
{\bf Remark 3.1.}
For more general condition on $D$ instead of it's convexity see \cite{IW3}. In that case,
instead of ($\ref{eq:3.3}$) we have to assume that
$$\displaystyle
\text{det}\,(S_{q}(E_c(p,p'))-S_{q}(\partial D))>0\,\,\forall q\in\Lambda_{\partial D}(p,p').
$$
Note that in that case
$\Lambda_{\partial D}(p,p')$ does not necessary consists of a single point and thus
its finiteness should be assumed.  As a result the right-hand side on ($\ref{eq:3.4}$) should be changed.

As a consequence of Theorem 3.2 we obtain two procedures for extracting
the curvatures and principle directions.

First we describe a procedure for extracting the curvatures
of $\partial D$ at $q(p,p')$ provided $q=q(p,p')$ is known.
Formula ($\ref{eq:3.4}$) gives the information
$\text{det}\,(S_{q}(E_c(p,p'))-S_{q}(\partial D))$ from $u_f$ on $B'\times\,]0,\,T[$.
Then by restricting $u_f$ on $B''\times\,]0,\,T[$ we obtain also from ($\ref{eq:3.4}$)
$\text{det}\,(S_{q}(E_{c-s}(p,p'+s\omega))-S_{q}(\partial D))$ for $\omega=(q-p')/\vert q-p'\vert$.
Using an analogous equation to ($\ref{eq:2.23}$) (Lemma 5.1 in \cite{IW3}), we see that
these two quantities construct a linear system with the Gauss curvature $K_{\partial D}(q)$
and a modification of the mean curvature of $\partial D$ at $q$, that is
$$
\displaystyle
\tilde{H}_{\partial D}(q;p,p')
\equiv H_{\partial D}(q)-\frac{S_{q}(\partial D)(\mbox{\boldmath $A$}_q(p)\times\mbox{\boldmath $A$}_q(p'))
\cdot
(\mbox{\boldmath $A$}_q(p)\times\mbox{\boldmath $A$}_q(p'))}{2(1+\mbox{\boldmath $A$}_q(p)\cdot\mbox{\boldmath $A$}_q(p'))},
$$
where $\displaystyle\mbox{\boldmath $A$}_q(x)=(q-x)/\vert q-x\vert$.
The system corresponds to ($\ref{eq:2.24}$) and always uniquely solvable.
Thus, by solving the system we obtain those two curvatures.
Therefore, one can obtain an approximate shape of $\partial D$ around $q(p,p')$.

Second we show that it is posible to obtain the principle curvature directions of
$\partial D$ at $q(p,p')$ by making a rotation of $B$ and $B'$ at the same time around the normal at $q(p,p')$.
We denote by $p(\theta)$ and $p'(\theta)$ the points rotated around the line directed $\nu_q$ at
$q=q(p,p')$ counterclockwise
with rotation angle $\theta\in[0,\,2\pi[$ of $p$ and $p'$.
Thus $p(0)=p$ an $p'(0)=p'$.
Then, for all $\theta\in\,[0,\,2\pi[$ we know that $\Lambda_{\partial D}(p(\theta),p'(\theta))=\{q(p,p')\}$,
$\mbox{\boldmath $A$}_q(p(\theta))\cdot\mbox{\boldmath $A$}_q(p'(\theta))$,
$\vert\mbox{\boldmath $A$}_q(p(\theta))\times\mbox{\boldmath $A$}_q(p'(\theta))\vert$
and $\vert p(\theta)-q\vert+\vert q-p'(\theta)\vert$ at $q=q(p,p')$
are invariant with respect to $\theta$.

Let $B(\theta)=\{x\in\Bbb R^3\,\vert\,\vert x-p(\theta)\vert<\eta\}$ and $B(0)=B$;
$B'(\theta)=\{x\in\Bbb R^3\,\vert\,\vert x-p'(\theta)\vert<\eta'\}$ and $B'(0)=B'$.
We have $[\overline B(\theta)\cup\overline B'(\theta)]\cap\overline D=\emptyset$
provided $[\overline B\cup\overline B']\cap\overline D=\emptyset$.
Let $f(\theta)$ denote the characteristic function of $B(\theta)$.

Then, from $u_{f(\theta)}$ on $B'(\theta)\times\,]0,\,T[$
we obtain the function of $\theta$:
$$\begin{array}{c}
\displaystyle
\tilde{H}_{\partial D}(q;p(\theta),p'(\theta))
=
H_{\partial D}(q)
-\frac{1-\mbox{\boldmath $A$}_q(p)\cdot\mbox{\boldmath $A$}_q(p')}{2}
S_{q}(\partial D)(\mbox{\boldmath $V$}(\theta))
\cdot
\mbox{\boldmath $V$}(\theta)
\end{array}
$$
where $\mbox{\boldmath $V$}(\theta)$ denotes the unit vector directed to
$\mbox{\boldmath $A$}_q(p(\theta))\times\mbox{\boldmath $A$}_q(p'(\theta))$.

Now assume that $\mbox{\boldmath $A$}_q(p)\times\mbox{\boldmath $A$}_q(p')\not=\mbox{\boldmath $0$}$.
Then, $\mbox{\boldmath $V$}(\theta)$ attains all the tangent vector at
$q$ of $\partial D$ and thus from the behaviour of $\tilde{H}_{\partial D}(q;p(\theta),p'(\theta))$
as a function of $\theta$ one can determine all the directions of principle curvatures say,
$\mbox{\boldmath $V$}(\theta_1)$ and $\mbox{\boldmath $V$}(\theta_2)$ with some $\theta_1$ and $\theta_2$.
Since we have
$$\displaystyle
S_{q}(\partial D)(\mbox{\boldmath $V$}(\theta_1))
\cdot
\mbox{\boldmath $V$}(\theta_1)+S_{q}(\partial D)(\mbox{\boldmath $V$}(\theta_2))
\cdot
\mbox{\boldmath $V$}(\theta_2)=2H_{\partial D}(q),
$$
the arithmetic mean of $\tilde{H}_{\partial D}(q;p(\theta_1),p'(\theta_1))$ and
$\tilde{H}_{\partial D}(q;p(\theta_2),p'(\theta_2))$
coincides with
$$\begin{array}{c}
\displaystyle
\left\{1-
\frac{1-\mbox{\boldmath $A$}_q(p)\cdot\mbox{\boldmath $A$}_q(p')}{2}\right\}
H_{\partial D}(q).
\end{array}
$$
Thus we obtain $H_{\partial D}(q)$.
Therefore, we can extract $S_q(\partial D)$
from $u_{f(\theta)}$ over $B'(\theta)\times\,]0,\,T[$ given at all $\theta\in [0,\,2\pi[$.
This is an advantage of the data collection using a {\it simultaneous rotation} of the emitter and the receiver.

\noindent
{\bf Open problem 3.1.} Extend the results to other boundary conditions, transmission conditions
(see Section 2.7) or the Maxwell system
(see Section 3.2).

\subsection{The Maxwell system}

In this section we briefly comment on a recent application \cite{IMAX} of the Enclosure Method to an inverse obstacle problem
whose govering equation is given by the Maxwell system in the time domain.

Let $0<T<\infty$.  We denote by $\mbox{\boldmath $E$}$
and $\mbox{\boldmath $H$}$ the elctric and magnetic fields, respectively.
Assume that $\mbox{\boldmath $E$}$ and $\mbox{\boldmath $H$}$
are induced only by the current density $\mbox{\boldmath $J$}$ at $t=0$ and that the obstacle is a perfect conductor
placed in the whole space $\Bbb R^3$.
The governing equations of $\mbox{\boldmath $E$}$ and $\mbox{\boldmath $H$}$
take the form
%$$\left\{
%\begin{array}{c}
%\displaystyle
%\epsilon\frac{\partial\mbox{\boldmath $E$}}{\partial t}
%-\nabla\times\mbox{\boldmath $H$}=\mbox{\boldmath $J$}
%\,\,\text{in}\,(\Bbb R^3\setminus\overline D)\times\,]0,\,T[,\\
%\\
%\displaystyle
%\mu\frac{\partial\mbox{\boldmath $H$}}{\partial t}
%+\nabla\times\mbox{\boldmath $E$}=\mbox{\boldmath $0$}\,\,\text{in}\,(\Bbb R^3\setminus\overline D)\times\,]0,\,T[,\\
%\\
%\displaystyle
%\mbox{\boldmath $\nu$}\times\mbox{\boldmath $E$}=\mbox{\boldmath $0$}\,\,\text{on}\,\partial D\times\,]0,\,T[,\\
%\\
%\displaystyle
%\mbox{\boldmath $E$}\vert_{t=0}=\mbox{\boldmath $0$},\,\,\mbox{\boldmath $H$}\vert_{t=0}=\mbox{\boldmath $0$}\,\,\text{in}\,\Bbb R^3\setminus\overline D,
%\end{array}
%\right.
%$$
$$
\left\{
\begin{array}{ll}
\displaystyle
\epsilon\frac{\partial\mbox{\boldmath $E$}}{\partial t}
-\nabla\times\mbox{\boldmath $H$}=\mbox{\boldmath $J$}
& \text{in}\quad (\Bbb R^3\setminus\overline D)\times\,]0,\,T[,
\\[10pt]
\displaystyle
\mu\frac{\partial\mbox{\boldmath $H$}}{\partial t}
+\nabla\times\mbox{\boldmath $E$}=\mbox{\boldmath $0$} & \text{in}\quad (\Bbb R^3\setminus\overline D)\times\,]0,\,T[,
\\[10pt]
\displaystyle
\mbox{\boldmath $\nu$}\times\mbox{\boldmath $E$}=\mbox{\boldmath $0$} & \text{on}\quad\partial D\times\,]0,\,T[,
\\[10pt]
\displaystyle
\mbox{\boldmath $E$}\vert_{t=0}=\mbox{\boldmath $0$} & \text{in}\quad\Bbb R^3\setminus\overline D,
\\[10pt]
\displaystyle
\mbox{\boldmath $H$}\vert_{t=0}=\mbox{\boldmath $0$} & \text{in}\quad\Bbb R^3\setminus\overline D,
\end{array}
\right.
$$
where $\mbox{\boldmath $\nu$}$ denotes the unit outward normal to $D$ on $\partial D$; $\epsilon$ and $\mu$ denote the electric permittivity
and magnetic permeability assumed to be positive constants.

There are several choices of the current density $\mbox{\boldmath $J$}$
as a model of the antenna.
In \cite{IMAX} it is assumed that $\mbox{\boldmath $J$}$
takes the form
%$$\displaystyle
%\mbox{\boldmath $J$}(x,t)
%=f(t)\chi_B(x)\mbox{\boldmath $a$},
%\tag {3.5}
%$$
\begin{eqnarray}\label{eq:3.5}
\displaystyle
\mbox{\boldmath $J$}(x,t)
=f(t)\chi_B(x)\mbox{\boldmath $a$},
\end{eqnarray}
where $\mbox{\boldmath $a$}\not=0$ is a constant unit vector,
$\chi_B$ denote the characteristic function of $B$ and $f\in
H^1(0,\,T)$ with $f(0)=0$; $B$ is an open ball with {\it very small} radius
and satisfies $\overline B\cap\overline D=\emptyset$.

In \cite{IMAX} the author considered the following problem.

\noindent
{\bf Problem 3.2.}  Fix $T$.  Generate $\mbox{\boldmath $E$}$ and $\mbox{\boldmath $H$}$ by the source $\mbox{\boldmath $J$}$ given by ($\ref{eq:3.5}$)
and observe $\mbox{\boldmath $E$}$ on $B$ over time interval $]0,\,T[$.
Extract information about the geometry of $D$ from the observed data.

Two theorems corresponding to Theorems 2.1 and 2.2 have been obtained in \cite{IMAX}.
The main difference from the scalar case is the existence of {\it directivity} of the source at $t=0$
and one of two theorems catches the effect of the source directivity.

The boundary condition imposed on the surface of the obstacle is a typical one like the Dirichlet boundary
condition for the wave equation.  As a next step
it is natural to ask: how about the case when the electromagnetic wave satisfies a more general
boundary condition like the Leontovich condition on the surface of the obstacle (see , e.g., \cite{ALN})?

\noindent
{\bf Open problem 3.2.}
Consider the Leontovich boundary condition instead of the perfect conductivity condition:
$$
\displaystyle
\mbox{\boldmath $\nu$}\times\mbox{\boldmath $H$}
-\lambda(x)\,
\mbox{\boldmath $\nu$}\times(\mbox{\boldmath $E$}\times\mbox{\boldmath $\nu$})=0\,\,\text{on}\,\partial D\times\,]0,\,T[.
$$
Extract information about the geometry of $D$ and $\lambda$ from the observed data.

$\quad$

\centerline{{\bf Acknowledgement}}

The author was partially supported by Grant-in-Aid for
Scientific Research (C)(No. 2540015) of Japan  Society for the
Promotion of Science.

%\vskip1cm
%\noindent
%e-mail address

%ikehata@amath.hiroshima-u.ac.jp


\begin{thebibliography}{99}
\setlength{\itemsep}{6pt}
\addtolength{\leftmargin}{-0.5cm}
\bibliographystyle{plain}
\footnotesize{




\bibitem{ALN} Ammari, H., Latiri-Grouz, C. and N\'ed\'elec, J. -C.,
              The Leontovich boundary value problem for the time-harmonic Maxwell equations,
              Asymptotic Analysis, {\bf 18}(1998), 33-47.





\bibitem{BH}  Bleistein, N. and Handelsman, R. A., Asymptotic expansions of integrals (New York: Dover Publications), 1986.




\bibitem{BKS} Brander, T., Kar, M. and Salo, M.,
              Enclosure method for the $p$-Laplace equation,
              Inverse Problems, {\bf 31}(2015) 045001.




\bibitem{CH} Courant, R. and Hilbert, D., Methoden der Mathematischen Physik, Vol. {\bf 2} (Berlin: Springer), 1937.







\bibitem{DL}  Dautray, R. and Lions, J-L., Mathematical analysis and numerical methods for
              sciences and technology,
              Evolution problems I, Vol. {\bf 5} (Berlin: Springer), 1992.





\bibitem{GT}  Gilbarg, D. and Trudinger, N. S., Elliptic partial differential equations
              of second order, second.ed. (Berlin:Springer), 1983.




\bibitem{G}  Grisvard, P., Elliptic Problems in Nonsmooth Domains (Boston: Pitman), 1985.






\bibitem{IP} Ikehata, M., Reconstruction of the shape of the inclusion by boundary measurements,
             Commun. in partial differential equations, {\bf 23}(1998), 1459-1474.




\bibitem{IPR} Ikehata, M., Reconstruction of an obstacle from the scattering amplitude at a fixed frequency,
              Inverse Problems, {\bf 14}(1998), 949-954.







\bibitem{IE3} Ikehata, M., Enclosing a polygonal cavity in a two-dimensional bounded domain from Cauchy data,
            Inverse Problems, {\bf 15}(1999), 1231-1241.


\bibitem{IE2} Ikehata, M., How to draw a picture of an unknown inclusion from boundary measurements.
            Two mathematical inversion algorithms,
            J. Inv. Ill-Posed Problems, {\bf 7}(1999), 255-271.




\bibitem{IP2} Ikehata, M., Reconstruction of obstacle from boundary measurements,
              Wave Motion, {\bf 30}(1999), 205-223.



\bibitem{IE} Ikehata, M., Reconstruction of the support function for inclusion from boundary measurements,
            J. Inv. Ill-Posed Problems, {\bf 8}(2000), 367-378.







\bibitem{IHEAT} Ikehata, M., Extracting discontinuity in a heat conductive body.
                One-space dimensional case, Appl. Anal., {\bf 86}(2007), no. 8, 963-1005.






\bibitem{ICUBO} Ikehata, M.,
                A remark on the enclosure method for a body with an unknown homogeneous background
                conductivity,  CUBO A Mathematical Journal, {\bf 10}(2008), no. 2, 31-45.







\bibitem{IW} Ikehata, M., The enclosure method for inverse obstacle scattering problems with dynamical
             data over a finite time interval, Inverse Problems, {\bf 26}(2010) 055010(20pp).











\bibitem{ISURVEY}  Ikehata, M., The probe and enclosure methods for inverse obstacle scattering problems.
                   The past and present., in {\it New Development of Functional Equations in Mathematical Analysis},
                   RIMS Kokyuroku, No. {\bf 1702}, 1-22, 2010. arXiv:1002.0198v4[math.AP].











\bibitem{IFRAME}  Ikehata, M., The framework of the enclosure method with dynamical data and its applications,
                  Inverse Problems, {\bf 27}(2011) 065005(16pp).







\bibitem{ILOG}  Ikehata, M., Inverse obstacle scattering problems with a single incident wave and the logarithmic differential of the
                indicator function in the enclosure method, Inverse Problems, {\bf 27}(2011) 085006(23pp).










\bibitem{IW2} Ikehata, M., The enclosure method for inverse obstacle scattering problems with dynamical data over a finite
              time interval: II. Obstacles with a dissipative boundary or finite refractive index and back-scattering
              data, Inverse Problems, {\bf 28}(2012) 045010(29pp).






\bibitem{IWcavity} Ikehata, M., An inverse acoustic scattering problem inside a cavity with dynamical back-scattering data,
                   Inverse Problems, {\bf 28}(2012) 095016(24pp).






\bibitem{Isugaku} Ikehata, M., Analytical methods for extracting discontinuity in inverse problems: the
                  probe method after 10 years, Sugaku Expositions, {\bf 26}(2013), Number 1, 1-28, AMS.







\bibitem{IW3} Ikehata, M., The enclosure method for inverse obstacle scattering problems with dynamical data over a finite time
              interval:  III. Sound-soft obstacle and bistatic data, Inverse Problems, {\bf 29}(2013) 085013(35pp).



\bibitem{IWCOE} Ikehata, M., Extracting the geometry of an obstacle and a zeroth-order coefficient of a boundary condition via the
                enclosure method using a single reflected wave over a finite time interval,
                Inverse Problems, {\bf 30}(2014) 045011(24pp).





\bibitem{IMAX} Ikehata, M., The enclosure method for inverse obstacle scattering
               using a single electromagnetic wave in time domain, Inverse Problems and Imaging, to appear.






\bibitem{IkIt} Ikehata, M. and Itou, H.,
               Extracting the support function of a cavity in an isotropic elastic body from a single set of boundary data,
               Inverse Problems, {\bf 25}(2009) 105005(21pp).



\bibitem{IkIt2} Ikehata, M. and Itou, H.,
               On reconstruction of a cavity in a linearized viscoelastic body from infinitely many transient boundary data,
               Inverse Problems, {\bf 28}(2012) 125003(19pp).




\bibitem{IK1}  Ikehata, M. and Kawashita, M., The enclosure method for the heat equation,
               Inverse Problems, {\bf 25}(2009) 075005(10pp).




\bibitem{IK2}  Ikehata, M. and Kawashita, M., On the reconstruction of inclusions in a heat conductive body
               from dynamical boundary data over a finite time interval, Inverse Problems,
               {\bf 26}(2010) 095004(15pp).





\bibitem{IK3}  Ikehata, M. and Kawashita, M., An inverse problem for a three-dimensional heat equation in
               thermal imaging and the enclosure method, Inverse Problems and Imaging, {\bf 8}(2014), No.4,
               1073-1116.










\bibitem{IO} Ikehata, M. and Ohe, T.,
             The enclosure method for an inverse crack problem and the Mittag-Leffler function,
             Inverse Problems, {\bf 24}(2008)015006(27pp).






\bibitem{IS} Ikehata, M. and Siltanen, S.,
              \newblock Electrical impedance tomography and Mittag-Leffler's function,
              \newblock Inverse Problems, 20(2004), 1325-1348.













\bibitem{ISA} Isakov, V., Inverse obstacle problems,
               Topical review, Inverse Problems, {\bf 25}(2009) 123002(18p).








\bibitem{LP}  Lax, P. D. and Phillips, R. S., The scattering of sound waves by an obstacle,
              Comm. Pure and Appl. Math., {\bf 30}(1977), 195-233.






\bibitem{SY} Sini, M. and Yoshida, K.,
             On the reconstruction of interfaces using complex geometrical optics solutions for the acoustic case,
             Inverse Problems, {\bf 28}(2012) 055013(22pp).




\bibitem{SY2} Sini, M. and Yoshida, K.,
              Corrigendum: On the reconstruction of interfaces using complex geometrical optics solutions for the acoustic case,
              2012 {\it Inverse Problems} {\bf 28} 055013,
              Inverse Problems, {\bf 29}(2013) 039501 (2pp).








\bibitem{WZ} Wang, J.-N. and Zhou, T.,
             Enclosure methods for Helmholtz-type equations.
             In {\it Inverse problems and applications: inside out. II},
             volume {\bf 60} of Math. Sci. Res. Inst. Publ., 249-270
             (Cambridge: Cambridge Univ. Press), 2013.





}















\end{thebibliography}
\end{document}